\documentclass[11pt,a4paper]{article}
\usepackage{a4wide,amsmath,amsthm,epsfig,graphicx}
\usepackage{amsmath,amsthm,amssymb}
\usepackage{amsfonts}
\usepackage{graphicx}
\title{\bf{Monte Carlo Methods and Path-Generation techniques for Pricing Multi-asset Path-dependent Options}}
\author{Piergiacomo Sabino\\
Dipartimento di Matematica\\
Universit\`{a} degli Studi di Bari\\
sabino@dm.uniba.it\\
Report 36/07}
\date{}
\begin{document}
\maketitle \thispagestyle{empty}
\abstract \noindent We consider the problem of pricing
path-dependent options on a basket of underlying assets using
simulations. As an example we develop our studies using Asian
options.

Asian options are derivative contracts in which the underlying
variable is the average price of given assets sampled over a period
of time. Due to this structure, Asian options display a lower
volatility and  are therefore cheaper than their standard European
counterparts.

This paper is a survey of some recent enhancements to improve
efficiency when pricing Asian options by Monte Carlo simulation in
the Black-Scholes model. We analyze the dynamics with constant and
time-dependent volatilities of the underlying asset returns.

We present a comparison between the precision of the standard Monte
Carlo method (MC) and the stratified Latin Hypercube Sampling (LHS).
In particular, we discuss the use of low-discrepancy sequences, also
known as Quasi-Monte Carlo method (QMC), and a randomized version of
these sequences, known as Randomized Quasi Monte Carlo (RQMC). The
latter has proven to be a useful variance reduction technique for
both problems of up to $20$ dimensions and for very high dimensions.

Moreover, we present and test a new path generation approach based
on a Kronecker product approximation (KPA) in the case of
time-dependent volatilities. KPA proves to be a fast generation
technique and reduces the computational cost of the simulation
procedure.
Path-dependent options.
\endabstract
%
\section{Introduction}
The financial industry has developed a variety of derivative
contracts in order to fulfil different investor needs.
Path-dependent options play a fundamental role in  financial
engineering  and can display different exotic features.

Exotic contracts that are widely used are Asian options, barrier
options and look-back options both with American and European style.
An un-biased and efficient pricing procedure is fundamental and a
vast research is done in order to obtain fast and efficient
estimations. Common approaches rely on finite differences methods
and Monte Carlo simulations.

Finite differences methods consist in discretizing the partial
differential equation whose solution gives the price of the options
while Monte Carlo methods face the problem from a probabilistic
point of view. It estimates the price as an expected value by its
integral formulation.

The former method returns the fair price of the option for different
times and values of the underlying variable but is practically
unfeasible for complicated multi-asset dependence.

On the other hand, Monte Carlo simulation calculates the fair price
in a single time point and can be applied to various situations.

Its fundamental property is that its order of convergence is
$O(1/\sqrt{n})$ and does not depend on the number of random sources
of the problem. Although it does not display a high order of
convergence, it proves to be efficient for pricing complex exotic
contracts.

The aim of this report is to describe standard and advanced Monte
Carlo techniques applied for multi-asset Asian options of European
style. In particular we concentrate our studies to stratification
and Quasi-Monte Carlo approaches.

Standard Monte Carlo can be seen as a numerical procedure aimed to
estimate integrals in the hypercube $[0,1]^d$ by generating
different scenarios with uniform random variables. Stratification
achieves the same task by drawing uniform random variates in a
smaller set in $[0,1]^d$ introducing correlation.

Quasi-Monte Carlo methods drop off all probabilistic considerations
and focus on the problem of generating a sequence of points that
uniformly covers the hypercube $[0,1)^d$ (the theory is built-up for
right-opened intervals). The sequence is absolutely deterministic
and different drawings lead to the same points.

From the mathematical point of view, it introduces the concept of
discrepancy and star-discrepancy that quantify how well the
sequences cover $[0,1)^d$.

Hawkla and Koksma proved the fundamental inequality, named after
them, that provides the bound for the estimation error of the target
integral depending on the discrepancy.

Low-discrepancy sequences are those whose estimation error is
$O(\frac{ln^d n}{n})$. The convergence rate depends on the dimension
$d$ and is lower than the Monte Carlo error for small $d$. There
exist several low discrepancy sequences, among them are the Halton,
the Faure, the Sobol and the Niederreiter-Xing sequences. A
fundamental reference on this topic is Niederreiter \cite{Ni1992}.

Quasi-Monte Carlo methods can be unpractical because the computation
of the error is potentially more difficult than the estimation of
the target integral while the better uniformity can be lost even for
low values of $d$.

Standard Monte Carlo, stratification and Quasi-Monte Carlo methods
form a hierarchy for the generation of uniform points.

A further step ahead can be taken by randomizing these sequences
while preserving the low-discrepancy. This technique is called
scrambling, Owen \cite{ow2002} provides an extensive description on
the subject.

The application to options pricing is straightforward. Standard
models for price dynamics involve multidimensional It\^{o} processes
so that pricing exotic contracts might require a high-dimensional
integration. It necessitates careful implementation of the
simulation especially when Quasi-Monte Carlo methods are used.

Many works have been done to investigate the problem. Acworth,
Broadie, and Glasserman \cite{ABG1998} provided  a first comparisons
between variance reduction techniques and Quasi-Monte Carlo methods
and Caflisch, Morokoff and Owen \cite{CMO1997} analyzed the
effective dimension of the integration problem for mortgage-backed
securities by ANOVA considerations. Caflisch, Morokoff and Owen
\cite{CMO1997} and Owen \cite{ow1998} showed that only few random
sources really matter and suggested to choose for them a better
generation technique.

We focus our investigation on pricing Asian options in a
multi-dimensional Black-Scholes model both for constant and
time-dependent volatilities. In this framework, standard
path-generating techniques are the Cholesky decomposition, the
principal component analysis (PCA)  and linear transform LT. The
last two have been proven to be essential for ANOVA in order to
recognize effective dimensions so that an efficient RQMC can be run.

When constant volatilities are considered, the path-generation
procedure can be simplified relying on the properties on the
Kronecker product while this is not possible for time-depending
volatilities.

As for this task, we propose a new approach based on a Kronecker
product approximation. The general problem consists in approximating
the global correlation matrix of the price returns into the
Kronecker product of two smaller matrices. We assume that the former
of the two is the auto-covariance matrix of a single brownian
motion. Indeed, we suppose that most of the variance of the global
process is carried out by each driving brownian motion. The latter
matrix would be an approximation of the total covariance matrix
among the asset returns during the lifetime of the contract. The
original and target path is reobtained by Cholesky decomposition. As
for this last step we develop an \emph{ad hoc} realization of the
Cholesky decomposition suited for the global correlation matrix.
This procedure is intended to reduce the computational burden
required to evaluate the whole set of eigenvalues and eigenvectors
of the global covariance matrix.

The last step of the simulation is the computation of the Asian
price via simulation using standard Monte Carlo, LHS and RQMC
approaches. As for the last one we perform a Faure-Tezuka scrabling
version of the Sobol' sequence, which is the most used low
discrepancy sequence in finance.

In the case of constant volatility, we set our investigation as in
Dahl and Benth \cite{DB2001}, \cite{DB2002} and Imai and Tan
\cite{IT2002}. We compare our results and analyze the precision of
the simulation for different path-generation methods and Monte
Carlo approach.

As for the time-dependent volatility market we test the KPA method
and compare its results with those obtained with the PCA
decomposition.

Summarizing the report evolves as follows: section 2 introduces the
market. Section 3 describes the pay-off of Asian options and
presents the problem as an integral formulation. Section 4 defines
effective dimensions in truncation and superposition sense. Section
5 defines the Kronecker product and list some of the main
properties. Section 6 describes some path-generation procedures and
in particular, it introduces the KPA  method. Section 7 is a brief
introduction to low-discrepancy sequences and scrambling techniques.
Section 8 describes the simulation procedure we adopt. Section 9
shows and comments the estimated results for different scenarios
both in the constant and time-dependent cases.
%
\section{The Market}
We consider a complete, standard financial market $\mathfrak{M}$ in
a Black-Scholes framework, with constant risk-free rate $r$ and
time-dependent volatilities. There are $M + 1$ assets in the market,
one risk free asset and $M$ risky assets. The price processes of the
assets in this market are driven by a set of stochastic differential
equations.

Suppose we have already applied the Girsanov
theorem and found the (unique) risk-neutral probability, the model
for the risky assets is the so called multi-dimensional geometric
brownian motion:
\begin{equation}
S_{0}(t)=e^{rt}\label{1.1}
\end{equation}
\begin{equation}
dS_{i}\left( t\right) =rS_{i}\left( t\right) dt+\sigma _{i}\left(
t\right)S_{i}\left( t\right)\,dW_{i}\left( t\right) ,\qquad
i=1,\dots ,M.  \label{1.2}
\end{equation}%
\noindent Here $S_{i}\left( t\right) $ denotes the $i$-th asset
price at time $t$, $\sigma _{i}\left( t\right)$ represents the
instantaneous time-dependent volatility of the $i$-th asset return,
$r$ is the continuously compounded risk-free  interest rate, and $%
\mathbf{W}\left( t\right) =\left( W_{1}\left( t\right) ,\dots
,W_{M}\left( t\right) \right) $ is an $M$-dimensional Brownian
motion. Time $t$ can vary in $\mathbb{R}_{+}^*$, that is, we can
consider any maturity $T\in\mathbb{R}_{+}^*$ for all financial
contracts.

The multi-dimensional brownian motion $\mathbf{W}\left( t\right)$ is
a martingale, that is, each component
is a martingale, and satisfies the following properties:%
\begin{equation*}
\mathbb{E}\left[ W_{i}\left( t\right) \right] =0,\qquad i=1,\dots
,M.
\end{equation*}%
\begin{equation*}
\left[W_{i},W_{k}\right]\left(t\right)=\rho _{ik}t,\qquad i,k =
1,\dots ,M.
\end{equation*}%
\noindent where $[\centerdot,\centerdot](t)$ represents the
quadratic variation up to time $t$ and $\rho _{ik}$ the constant
instantaneous correlation between $W_{i}$ and $W_{k}$.

Consider a generic maturity $T$, we can define a time grid
$\mathcal{T} = \{t_1,\dots,t_N\}$ of $N$ points such that
$t_1<t_2<\dots,t_N=T$, we recall that the sampled covariance matrix
$R_{l,m}= \mathbb{E}\left[ W_i\left( t_{l}\right) W_i\left(
t_{m}\right) \right]$, $l,m=1,\dots ,N$
of each Brownian motion in equation (\ref{1.2}) is:%
\begin{equation}\label{1.3}
R=\left(
\begin{array}{cccc}
t_{1} & t_{1} & \ldots & t_{1} \\
t_{1} & t_{2} & \ldots & t_{2} \\
\vdots & \vdots & \ddots & \vdots \\
t_{1} & t_{2} & \ldots & t_{N}%
\end{array}%
\right)\text{.}
\end{equation}%
\noindent This matrix is symmetric and its elements $R_{l,m}=
t_l\wedge t_m$ have the peculiarity to be constant after reflection
about the diagonal. We will refer to this feature as
\emph{boomerang} shape property.

In order to complete the picture of our environment, we need to
define the matrix $\Sigma(t)$, whose elements are ${\Sigma
}_{i,k}(t)$ =$\rho _{ik}\sigma _{i}(t)\sigma _{k}(t),\,i,k=1,\dots
\,M$. This is a time dependent covariance matrix evolving according
to the dynamics of the time-dependent volatilities and the constant
correlation among the asset returns.

Avoiding all the calculation (see  Rebonato \cite{Re2004} and
Glassermann \cite{Glass2004} for further details), we derive the
global covariance matrix $\Sigma_{MN}$ that assumes the expression
below:
\begin{equation}
\Sigma_{MN} = \left( \begin{array}{cccc}
\Sigma(t_1) & \Sigma(t_1) & \ldots &
\Sigma(t_1)
\\ \Sigma(t_{1}) & \Sigma(t_2)& \ldots & \Sigma(t_2)
\\ \vdots & \vdots & \ddots & \vdots
\\\Sigma(t_{1}) & \Sigma(t_2) & \ldots & \Sigma(t_{N})
\end{array} \right)\label{1.4}
\end{equation}
The global covariance matrix is very simple and enjoys the
\emph{boomerang} shape property with respect to the block-matrix
notation. All the information is carried out by $N$ time-varying
$M\times M$ matrices.

Each element depends on four indexes:
\begin{equation}
\Big(\big(\Sigma_{MN}\big)_{ik}\Big)_{lm} = \int_0^{t_l\wedge
t_m} \sigma_i(t) \sigma_k(t)\rho_{ik}dt\label{1.5}
\end{equation}
with $i,k=1,\dots,M$ and $l,m=1,\dots,N$.

Applying the risk-neutral pricing formula, the value at time $t$ of
any European $T$-maturing derivative contract is:
\begin{equation}
V(t) =
exp\left(r(T-t)\right)\mathbb{E}\left[\phi(\mathcal{T})\right|\mathcal{F}_t]\label{1.6}\text{.}
\end{equation}
\noindent $\mathbb{E}$ denotes the expectation under the risk
neutral probability measure and $\phi(\mathcal{T})$ is a generic
$\mathcal{F}_T$ measurable function, with
$\mathcal{F}_T=\sigma\{0<t\leq T;W(t)\}$, that determines the payoff
of the contract. Although not explicitly written, the function
$\phi(\mathcal{T})$ depends on the entire multi-dimensional brownian
path up to time $T$.
\section{Problem Settlement}
We will restrict our analysis to Asian options that  are exotic
derivative contracts that can be written both on a single security
and on a basket of underlying securities. Hereafter we will consider
European-style Asian options whose underlying securities coincide
with the $M+1$ assets on the market. This is the most general case
we can tackle in the market $\mathfrak{M}$, because it is complete
in the sense that we can hedge any financial instrument by finding a
portfolio that is a combination of this $M+1$ assets.
\subsection {Asian Options Payoff}
The theoretical definition of Asian options price is:
\begin{equation}
    a_i(t) = exp\left(r(T-t)\right)\mathbb{E}\left[\left(\frac{\int_0^TS_i(t)dt}{T}-K\right)^+\bigg|\mathcal{F}_t\right]
    \label{2.1.1}\quad\text{Option on a Single Asset}
\end{equation}
\begin{equation}
    a(t) =
    exp\left(r(T-t)\right)\mathbb{E}\left[\left(\frac{\int_0^T\sum_{i=1}^M
    w_iS_i(t)dt}{T}-K\right)^+\bigg|\mathcal{F}_t\right],
    \label{2.1.2}\quad\text{Option on a Basket}
\end{equation}
\noindent where we assume that the start date of the contract is
 $t=0$. $K$ represents the strike price and coefficients $w_{i}$
satisfy $\sum_{i=1}^Mw_{i}=1$. Contingent claims ($\ref{2.1.1}$) and ($\ref{2.1.2}$) are
usually referred as weighted Asian options.

In practice no contract is agreed according to equations
(\ref{2.1.1}) and (\ref{2.1.2}). The integrals are approximated by
sums; often these approximations are written in the contracts by
specifying the number and the sampling points of the path.

Approximation for (\ref{2.1.1}) and (\ref{2.1.2}) can be carried out
by using the following expressions:
\begin{equation}
a_i\left( t\right) =exp\left(r(T-t)\right)\mathbb{E}\left[\left(
\frac{\sum_{j=1}^{N}S_{i}\left( t_{j}\right)}{N}
-K\right)^+\bigg|\mathcal{F}_t\right] \quad\text{Option on a Single
Asset}\label{2.1.3}
\end{equation}
\begin{equation}
a\left( t\right) =exp\left(r(T-t)\right)\mathbb{E}\left[\left(
\sum_{i=1}^{M}\sum_{j=1}^{N}w_{ij}\,S_{i}\left( t_{j}\right)
-K\right)^+\bigg|\mathcal{F}_t\right]\quad\text{Option on a Basket}
\label{2.1.4}
\end{equation}%
\noindent where coefficients $w_{ij}$ satisfy $\sum_{i,j}w_{ij}=1$.

European options with payoff functions $\left( \ref{2.1.3}\right)$
and $\left( \ref{2.1.4}\right)$ are called arithmetic weighted
average options or simply arithmetic Asian options. When $M>0$ and
$N=1$ the payoff only depends on the terminal price of the basket of
$M$ underlying assets and the option is known as basket option.

No closed-form solution exists for Asian options arbitrage-free
price, neither for single nor for basket options both for
theoretical and finitely monitored payoff. In order to obtain a
correct valuation of the price we are compelled to turn to numerical
procedures such as the Monte Carlo estimation or the finite
difference methods.

The latter is based on a convenient and correct discretization of
the partial differential equation associated to the risk neutral
pricing formula via the Feynmann-Kac representation. The finite
difference method returns the price for all the times and initial
values of the underlying assets. Ve\^{c}er \cite{Vecer2001} and
\cite{Vecer2002} found a convenient approach for the single asset
case and presents the comparison with other techniques. The main
drawback is the stability of the method that is practically
unfeasible for options on a basket.

Monte Carlo simulation is a numerically intensive methodology that
provides unbiased estimates with convergence rate not depending on
the dimension of the problem (the number of random sources to draw).
The cases of high values for the problem dimension  find interesting
applications in finance including the pricing of high-dimensional
multi-factor path-dependent options. In contrast to the finite
difference technique, the Monte Carlo method returns the estimate
for a single point in time. It is a flexible approach but requires
\emph{ad hoc} implementation and refinements, such as variance
reduction techniques, in order to improve its efficiency.

The main purpose of the standard Monte Carlo method is to
numerically estimate the integral below:
\begin{equation}
I = \int_{[0,1]^d} f(\mathbf{x})\,\mathbf{dx} \label{2.1.5}\text{.}
\end{equation}
\noindent The integral $I$ can be regarded as
$\mathbb{E}\left[f(U)\right]$, the expected value of a function $f(\centerdot)$ of the random
vector $\mathbf{U}$ that is uniformly distributed in hypercube $[0,1]^d$.

Monte Carlo methods simply estimate $I$ by drawing a sample of $n$
independent replicates $U_1\dots,U_n$ of $\mathbf{U}$ and then
computing the arithmetic average:
\begin{equation}
\widehat{I} = \widehat{I}_n = \frac{1}{n} \sum_{i=1}^n
f(U_i).\label{2.1.6}
\end{equation}

The Law of Large Numbers ensures that $\widehat{I}_n$ converges to
$I$ in probability a.s. and the Central Limit Theorem states that $I
- \widehat{I}_n$ converges in distribution to a normal with mean $0$
and standard deviation $\sigma/\sqrt{n}$ with
$\sigma=\sqrt{\int_0^1\left(f(\mathbf{x})-I\right)^2\mathbf{dx}}$.
The convergence rate is than $O(1/\sqrt{n})$ for all dimensions $d$.
The parameter $\sigma$ is generally unknown in a setting in which
$I$ is unknown, but it can be estimated using the sampled standard
deviation or root mean square error (RMSE):
\begin{equation}
 RMSE = \sqrt{\frac{1}{n-1}\sum_{i=1}^n\left( f(U_i) -
 \widehat{I}_n\right)^2}.\label{2.1.7}
\end{equation}
Refinements in Monte Carlo methods consist in finding techniques
whose aim is to reduce the RMSE, known as variance reduction
techniques, without changing the convergence rate. In contrast, the
Quasi Monte Carlo version focuses on the improvement of the
convergence rate by generating sequences in $[0,1]^d$ with high
stratification in order to uniformly cover  the hypercube. These
sequences are no longer random  and estimates and errors are not
based on probabilistic considerations.

As far as our case is concerned, we need to formulate the problems
(\ref{2.1.3}) and (\ref{2.1.4}) for pricing Asian options as
integrals of the form (\ref{2.1.5}) in order to apply the Monte
Carlo procedure.
\subsection {Problem Formulation as an Integral}
The model $\mathfrak{M}$, presented in the first section, consists
of the risk-free money market account and $M$ assets driven $M$
geometric brownian motion described by equation (\ref{1.2}) whose
solution is:
\begin{equation}
S_{i}\left( t\right) =S_{i}\left( 0\right) exp\left[ \int_0^t\left(
r- \frac{\sigma _i^2\left( s\right)}{2}\right)ds +\int_0^t\sigma
_i\left( s\right)dW_{i}\left( s\right) \right] ,i=1,...,M.\label{2.2.1}
\end{equation}
\noindent The quantity $\int_0^T\frac{\sigma _i^2\left(
s\right)}{T}ds$ is the total volatility for the $i$-th asset. The
solution (\ref{2.2.1}) is a multi-dimensional geometric brownian
motion, written GBM$\left(r,\int_0^t\frac{\sigma _i^2\left(
s\right)}{2}ds\right)$, in the sense that it can be obtained
applying It\^{o}'s lemma to
$S_i(t)=f\left(X_i(t)\right)=e^{X_i(t)}$, with $X_i(t)$ the $i$-th
component of the multi-dimensional brownian motion with drift $r$
and $i$-th diffusion $\int_0^t\frac{\sigma _i^2\left(
s\right)}{2}ds$, written BM$\left(r,\int_0^t\frac{\sigma _i^2\left(
s\right)}{2}ds\right)$.

Under the assumption of constant volatility the solution is still a
multi-dimensional geometric brownian motion with the following form:
\begin{equation}
S_{i}\left( t\right) =S_{i}\left( 0\right) exp\left[ \left( r-
\frac{\sigma _i^2}{2}\right)t +\sigma _iW_{i}\left( t\right) \right]
,i=1,...,M.\label{2.2.2}
\end{equation}
\noindent In compacted notation the solution (\ref{2.2.2}) is
GBM$\left(r,\frac{\sigma _i^2}{2}\right)$.

Pricing Asian option requires to monitor the solutions (\ref{2.2.1})
and (\ref{2.2.2}) at a finite set of points in time
$\{t_1,\dots,t_N\}$. This sampling procedure yields the following
expressions for time-dependent and constant volatilities:
\begin{equation}
S_i(t_j) = S_i(0)exp\bigg[\int_0^{t_j}\left(r -
\frac{\sigma^2_i(s)}{2}\right)ds + Z_i(t_j)\bigg]\label{2.2.3}
\end{equation}
\begin{equation}
S_i(t_j) = S_i(0)exp\bigg[\left(r -
\frac{\sigma^2_i}{2}\right)t_{t_j} + Z_i(t_j)\bigg]\label{2.2.4}
\end{equation}
\noindent where the components of the vector $\left(Z_1(t_1),\dots
Z_1(t_N),Z_2(t_1),\dots,Z_M(t_N)\right)$ are  $M\times N$ normal
random variables with zero mean vector and covariance matrix
$\Sigma_{MN}$, whose form simplifies in the case of constant
volatilities as we will be shown in Section 4.

The payoff at maturity $T$ of the arithmetic average Asian option is
then:
\begin{equation}
p_a(T) = \left(g(\mathbf{Z})-K\right)^+\label{2.2.5}
\end{equation}
where
\begin{equation}
g(\mathbf{Z})=\sum_{k=1}^{M\times N} exp\left (\mu_k + Z_k
\right)\label{2.2.6}
\end{equation}
and
\begin{equation}
\mu_k = \ln(w_{k_1k_2}S_{k_1}(0)) + \bigg(r-\frac{\sigma_{k_1}^2}{2}
\bigg)t_{k_2}\label{2.2.7}
\end{equation}
for constant volatilities or
\begin{equation}
\mu_k = \ln(w_{k_1k_2}S_{k_1}(0)) + rt_{k_2}-\frac{\int_0^{t_{k_2}}
\sigma_{k_1}^2(s)ds}{2}\label{2.2.8}
\end{equation}
for  time-dependent volatilities. The indexes $k_1$ and $k_2$ are
 $k_1=(k-1)modM, k_2 = [(k-1)/M]+1$, respectively, where $mod$ denotes
the modulus and $[\centerdot]$  the greatest integer less than or
equal to $x$.

The calculation of the price $a(t)$ in equation ($\ref{2.1.4}$)
can be formulated as an integral on $\mathbb{R}^{NM}$ in the
following way (see Dahl and Benth \cite{DB2001} and
\cite{DB2002}):
\begin{equation}\label{2.2.9}
    a\left( t\right) =exp\left(r(T-t)\right)
    \int_{\mathbb{R}^{MN}}\left(g(\mathbf{z})-K\right)^+F_{\mathbf{Z}}(\mathbf{dz})
\end{equation}
$F_{\mathbf{Z}}$ is the cumulative distribution of the normal random
vector $N(0,\Sigma_{NM})$.

In the following section we will show how to obtain the random
vector $\mathbf{Z}$ starting from a vector of independent and
normally distributed random variables $\mathbf{\epsilon}$. Once this
generation is carried out, we can apply the inverse transform method
to formulate the pricing problem as an integral of uniform random
variables in the hypercube $[0,1]^{MN}$ and use Monte Carlo methods:
\begin{equation}
    a\left( t\right) =exp\left(r(T-t)\right)
    \int_{[0,1]^{NM}}\left(g(\mathbf{u})-K\right)^+F^{-1}_{\mathbf{Z}}(\mathbf{u})\mathbf{du}\label{2.2.10}
\end{equation}

 In the following sections we will
present recent enhancement based on ANOVA for high dimensional Monte
Carlo and Quasi-Monte Carlo simulations in order to estimate the
integral ($\ref{2.2.10}$) for the pricing Asian option on a basket
of underlying assets both for constant and time-dependent
volatilities.
\section{Effective Dimensions}
When the nominal dimension $d$ of the problem of estimating the
integral $ (\ref{2.1.5})$ is one, there are standard numerical
techniques that give a good accuracy when $f$ is smooth.
Considerable problems arise when $d$ is high.

Recent studies proved that many financial experiments present
problem dimensions lower than the nominal one.  Owen (1998)
\cite{ow1998} and Caflisch, Morokoff and Owen \cite{CMO1997}
studied the application of ANOVA for high-dimensional problems and
introduced the definition of effective dimension. It is possible
to study some mathematical properties of the function $f$ and try
to split it in order to reduce the computational effort. The ANOVA
decomposition consists of finding a representation of $f$ into
orthogonal functions each of them depending only on a subset of
the original variables. This is the peculiar and stronger
condition that makes ANOVA different and more powerful with
respect to the usual Least Squared method.

Let $\mathcal{A}=\{1,\dots,d\}$ denote the set of the independent
variables for $f$ on $[0,1]^d$. $f$ could be written into the sum of
$2^d$ orthogonal functions each of them defined in different subsets
of $\mathcal{A}$, that is depending only on the variables in each of
these subsets:
\begin{equation}
f(\mathbf{x}) = \sum_{u\subseteq \mathcal{A}} f_u(\mathbf{x})\label{3.1}
\end{equation}
Now  let $|u|$ denote the cardinality of $u$, $\mathbf{x_u}$ the
$|u|$-tuple consisting of components $x_j$ with $j\in u$, and $-u$
being the complement of $u$ in $\mathcal{A}$. Then set the function
as:
\begin{equation}
f_u(\mathbf{x}) = \int_{\mathbf{z}:z_u=x_u} \bigg(f(\mathbf{z}) - \sum_{v\subset u}
f_v(\mathbf{z})\bigg)\,\mathbf{dz_{-u}}\label{3.2}
\end{equation}
Equation (\ref{3.2}) defines $f_u$ by subtracting what can be
attributed to the subsets of $u$, and then averaging over all
components not in $u$. In the function setting  $f_u(\mathbf{x_u})$
only depends on $\mathbf{x_u}$.\\ Denoting $\sigma^2 =
\int(f(\mathbf{x})-I)^2\,\mathbf{dx}$, $\sigma^{2}_{u} = \int
f_u(\mathbf{x})^2\,\mathbf{dx}$, $\sigma^{2}_{0}=0$, supposing
$\sigma < +\infty$ and $|u|> 0$ it follows:
\begin{equation}
\sigma^2 = \sum_{u\subseteq \mathcal{A}} \sigma^{2}_{u}\label{3.3}
\end{equation}
Equation (\ref{3.3}) partitions the total variance into parts
corresponding to each subset $u\subseteq\mathcal{A}$. The $f_u$
exhibits some nice properties: if $j\in u$ the line integral
$\int_{[0,1]} f_u(\mathbf{x})\,dx_j = 0$ for any $x_k$ with $k\neq
j$, and if $u\neq v$  $\int f_u(x) f_v(x)\,dx = 0$.

Exploiting the ANOVA decomposition the definition of effective
dimension can be given in the following ways:
\newtheorem{Definition1}{Definition}
\begin{Definition1}
The effective dimension of $f$, in the superposition sense, is the
smallest integer $d_S$ such that $\sum_{0<|u|\leq d_S} \sigma^2_u
\geq p \sigma^2$ .\\ The value $d_S$ depends on the order in
which the input variables are indexed.
\end{Definition1}
\newtheorem{Definition2}[Definition1]{Definition}
\begin{Definition2}
The effective dimension of $f$, in the truncation sense, is the
smallest integer $d_T$ such that $\sum_{u\subseteq\{1,\dots,d_t\}}
\sigma^2_u \geq p \sigma^2$.
\end{Definition2}
\noindent $0<p<1$ is an arbitrary level; the usual choice is $p=99\%$.

The definition of effective dimension in truncation sense reflects
that for some integrands only a small number of the inputs might
really matter. The definition of effective dimension in
superposition sense takes into account that for some integrands the
inputs might influence the outcome through their joint action within
small groups. Direct computation leads: $d_S \leq d_T \leq d$.
%
%
\section{The Kronecker Product}
The Black-Scholes model was originally built up under the hypothesis
of constant volatilities for all the assets. If this assumption
drops off the main ideas underlying the market $\mathfrak{M}$
described above do not change and fundamental results still hold.
The constant volatility case reduces the computational complexity of
the analysis and simplifies many calculations.

In the following we present some useful properties of the brownian
motion, of its sampled autocovariance matrix and of the global
covariance matrix. Furthermore, we introduce the Kronecker product
that will prove to be a powerful tool for reducing the computation
burden and a fast way to generate multi-dimensional brownian paths.

The sampled covariance matrix of each brownian motion, $R$, enjoys
many properties due to its particular \emph{boomerang} form. We list
some of them below:
\begin{enumerate}
\item The inverse of $R$ is a symmetric tri-diagonal matrix:
\begin{equation}\label{4.1}
\begin{array}{c}
  R^{-1} = \\
  \left( \begin{array}{cccccc}
\frac{t_{2}}{t_1(t_2-t_1)} & -\frac{1}{t_2-t_1} & 0 & \ldots &
\ldots & 0
\\ -\frac{1}{t_2-t_1} & \frac{t_3-t_1}{(t_2-t_1)(t_3-t_2)} & -\frac{1}{t_3-t_2} & 0 &
\ldots & \vdots
\\ 0 & -\frac{1}{t_3-t_2} &
\frac{t_4-t_2}{(t_3-t_2)(t_4-t_3)} & -\frac{1}{t_4-t_3} & \ldots &
\vdots
\\ \vdots & 0 &  -\frac{1}{t_4-t_3} &
 \ddots & \ddots & \vdots
\\ \vdots & \vdots & \ddots & \ddots & \frac{t_n-t_{n-2}}{(t_{n-1}-t_{n-2})(t_n-t_{n-1})} &
-\frac{1}{t_n-t_{n-1}}
\\ 0 & 0 & 0 & 0 & -\frac{1}{t_n-t_{n-1}} & \frac{1}{t_n-t_{n-1}}
\end{array} \right)
\end{array}
\end{equation}
$R^{-1}$ is a sparse matrix and low memory is required to store it.
$R$ and $R^{-1}$ share the same set of eigenvectors and have inverse
eigenvalues (the matrices are both definite positive).
\item The Cholesky decomposition of $R$ gives a \emph{boomerang} shaped matrix $C$.

\newtheorem{Cholesky}[Definition1]{Definition}
\begin{Cholesky}[\textbf{Cholesky Decomposition}]
Given any hermitian, definite positive matrix $A$, then $A$ can be decomposed as:
\end{Cholesky}
\begin{equation}
A = C_A\,C_A^*\label{4.2}
\end{equation}
where $C_A$ is a lower triangular matrix with strictly positive diagonal entries,
and C* denotes the conjugate transpose of C. The Cholesky decomposition is unique and the Cholesky matrix
can be interpreted as a sort of square root of $R$;
as far as the Cholesky decomposition of a symmetric matrix $A$ is concerned $C_A^*$ must be replaced by $C_A^T$.

After direct computation $C_R$ shows the form below:
\begin{equation}
C_R = \left( \begin{array}{cccc} \sqrt{t_{1}} & 0 & \ldots &
0
\\ \vdots & \sqrt{t_2-t_1} & \ddots & 0
\\ \vdots & \vdots & \ddots & \vdots
\\ \sqrt{t_{1}} & \sqrt{t_2-t_1} & \ldots & \sqrt{t_{N}-t_{N-1}}
\end{array} \right)\label{4.3}
\end{equation}
In the case of an equally spaced time grid, the Cholesky matrix is
just a lower triangular matrix whose elements are all equal to the
time step $\Delta t$.

\item The inverse of the Cholesky matrix is a sparse matrix, in
particular it is a bi-diagonal matrix whose elements on the same
row are equal and in opposite sign:
\begin{equation}
C_R^{-1}= \left(\begin{array}{ccccc} \frac{1}{\sqrt{t_{1}}} &
0 & \ldots & \ldots & 0
\\ -\frac{1}{\sqrt{t_2-t_1}} & \frac{1}{\sqrt{t_2-t_1}} & 0 & \dots & 0
\\ 0 & -\frac{1}{\sqrt{t_3-t_2}} & \frac{1}{\sqrt{t_3-t_2}} &
\ddots & \vdots
\\ \vdots & \vdots & \ddots & \ddots & \vdots
\\ 0 & 0 & 0 & -\frac{1}{\sqrt{t_n-t_{n-1}}} & \frac{1}{\sqrt{t_n-t_{n-1}}}
\end{array} \right)\label{4.4}
\end{equation}
\end{enumerate}
All these results prove to be useful for the simulation and reduce
the number of operations for the brownian path generation.

As for constant volatilities, both the covariance matrix among the
asset returns and the global covariance matrix simplify and are not
time-depending anymore.

Let $\Sigma$ be a covariance matrix depending on the
correlation among the asset returns whose elements are:
$\Sigma_{i,k}$ =$
\rho_{ik}\sigma_i\sigma_k,\,i,k=1,\dots\,M$, then
the global covariance matrix $\Sigma_{MN}$ displays the following form:
\begin{equation}
\Sigma_{MN} = \left( \begin{array}{cccc} t_{1}\Sigma &
t_{1}\Sigma & \ldots & t_{1}\Sigma
\\ t_{1}\Sigma & t_2\Sigma & \ldots & t_{2}\Sigma
\\ \vdots & \vdots & \ddots & \vdots
\\t_{1}\Sigma & t_2\Sigma & \ldots & t_{N}\Sigma
\end{array} \right)\label{4.5}
\end{equation}
This matrix is obtained by repeating the constant block of covariance
$\Sigma$ at all the points of the time grid.

This kind of mathematical operation is known as Kronecker product, denoted as $\otimes $. As such, $%
\Sigma_{MN}$ can be identified as the Kronecker product between $%
R$ and $\Sigma$, $R\otimes \Sigma$.
The Kronecker product reduces the computational complexity by
enabling operations on a $\left( N\times M,N\times M\right)$ matrix using
two smaller matrices that are $N\times N,$ and $M\times M$ respectively.
\newtheorem{Kronecker}[Definition1]{Definition}
\begin{Kronecker}[\textbf{The Kronecker Product}]
The Kronecker product of $A_{m_{A}\times n_A} \in
\mathbb{R}^{m_A\times n_A}$ and $B_{m_B\times n_B}\in\mathbb{R}^{m_B\times n_B}$,
written $A\otimes B$, is the tensor algebraic operation defined as:
\end{Kronecker}
\begin{equation}
A\otimes B = \left( \begin{array}{cccc} a_{11}B & a_{12}B &
\ldots & a_{1n_A}B
\\ a_{21}B & a_{22}B & \ldots & a_{2n_A}B
\\ \vdots & \vdots & \ddots & \vdots
\\a_{m_{A}1}B & a_{m_{A}1}B & \ldots & a_{m_{A}n_A}B
\end{array} \right)\label{4.6}
\end{equation}
The Kronecker product offers many properties some of these listed
below (for further details and proofs see Golub and Van Loan
\cite{GV1996}, Van Loan \cite{VL2000}, A.N. Langville, W.J.
Stewart \cite{LS2002}):
\begin{enumerate}
\item{Associativity.}
\begin{equation*}
A\otimes (B\otimes C) = (A\otimes B)\otimes C)
\end{equation*}
\item{Distributivity.}
\begin{equation*}
(A+B)\otimes (C+D) = A \otimes C + B \otimes C + A \otimes D + B
\otimes D
\end{equation*}
\item{Compatibility with ordinary matrix multiplication.}
\begin{equation*}
AB\otimes CD = (A\otimes C)(B\otimes D)
\end{equation*}
\item{Compatibility with ordinary matrix inversion.}
\begin{equation*}
(A\otimes B)^{-1} = A^{-1} \otimes B^{-1}
\end{equation*}
\item{Compatibility with ordinary matrix transposition.}
\begin{equation*}
(A\otimes B)^{T} = A^{T} \otimes B^{T}
\end{equation*}
\item{Trace factorization}
\begin{equation*}
tr(A\otimes B) = tr(A) tr(B)
\end{equation*}
\item{Norm factorization}
\begin{equation*}
\|A\otimes B\| = \|(A)\| \|(B)\|
\end{equation*}
\item{Compatibility with Cholesky decomposition.}
\\ Let $A$ and $B$ semi-definite positive matrices then:
\begin{equation*}
A\otimes B = (C_{A}C_{A}^{T}) \otimes (C_{B}C_{B}^{T}) =
(C_A\otimes C_B)(C_A\otimes C_B)^T
\end{equation*}
\item{Special matrices.}
\\ Let $A$ and $B$ be nonsingular, lower (upper) triangular, banded,
symmetric, positive definite, \dots, etc, then $A\otimes B$ preserves the
property.
\item{Eigenvalue and Eigenvectors.}\\
Define two square matrices $A$ and $B$, $N\times N$ and  $M\times
M$, respectively. Suppose that $\lambda_1,\dots,\lambda_N \in \sigma
(A)$, $\mathbf{v_1},\dots,\mathbf{v_N}$ and $\mu_1,\dots,\mu_M \in
\sigma (B)$,  $\mathbf{w_1},\dots,\mathbf{w_M}$ are the eigenvalues
and the correspondent eigenvectors of the two matrices respectively,
where $\sigma(\centerdot)$ denotes the spectre of the matrix. The
Kronecker product, $A\otimes B$, has eigenvectors
$\mathbf{v_i}\otimes \mathbf{w_j}$ and eigenvalues $\lambda_i
\mu_j$.

Summarizing, every eigenvalue of $A\otimes B$ arises as product of
eigenvalues of $A$ and $B$, and every eigenvector as a Kronecker
product between the corresponding eigenvectors. This last property
still holds for singular value decomposition.
\end{enumerate}
%
%
\section{Generating Sample Path}
In discussing the simulation of a geometric brownian motion we should
focus on the realization of a simple brownian motion at the sample time points
of the grid.

Because brownian motion has independent and normally distributed
increments, simulating $W_i(t_l)$ is straightforward.

Let $\epsilon_1,\dots,\epsilon_N$ be independent standard normal random variables and set $W_i(t_0)=0$.
Subsequent values can be generated as follow :
\begin{equation}
    W_i(t_l) = W_i(t_{l-1}) + \sqrt{t_l-t_{l-1}}\epsilon_l, l = 1,\dots,N\label{5.1}
\end{equation}
\noindent For a brownian motion $X_i(t)=$BM$(\mu_i,\sigma_i)$ given $X_i(t_0)$ set
\begin{equation}
    X_i(t_l) = X_i(t_{l-1}) +\mu_i\left(t_l-t_{l-1}\right)+
     \sqrt{t_l-t_{l-1}}\sigma_i\epsilon_l, l = 1,\dots,N\label{5.2}
\end{equation}
\noindent For time-dependent parameters the recursion becomes (in
the general situation the drift can be time-dependent too):
\begin{equation}
    X_i(t_l) = X_i(t_{l-1}) +\int_{t_{l-1}}^{t_l}\mu_i \left(s\right)ds+
    \sqrt{\int_{t_{l-1}}^{t_l}\sigma_i^2\left(s\right)ds}\epsilon_l, l = 1,\dots,N\label{5.3}
\end{equation}
\noindent The methods (\ref{5.1})-(\ref{5.3}) are exact in the sense
that the joint distribution of the random vector
$\left(W_i(t_1),\dots,W_i(t_N)\right)$ or
$\left(X_i(t_1),\dots,X_i(t_N)\right)$ coincides with that of the
original process at the times $\{t_1,\dots,t_N\}$, but are subject
to a discretization error.

Nothing can be said about what happens between the time point of the grid.
One might choose a linear interpolation to get intermediate values of the simulated process
without obtaining a correct joint distribution.

Applying the Euler scheme for the brownian motion with
time-dependent drift and diffusion,
\begin{equation}
X_i(t_l) = X_i(t_{l-1}) +\mu_i(t_l)\left(t_l-t_{l-1}\right)+
\sqrt{t_l-t_{l-1}}\sigma_i(t_l)\epsilon_l, l = 1,\dots,N\label{5.4}
\end{equation}
\noindent we introduce a dicretization error even at time points
$\{t_1,\dots,t_N\}$, because the increments will no longer have  the
correct mean and variance.

The vector $\left(W_i(t_1),\dots,W_i(t_N)\right)$ is a linear
combination of the vector of the increments
$\left(W_i(t_1)-W_i(t_0),\dots,W_i(t_N)-W_i(t_{N-1})\right)$ that is
normally distributed. All linear combinations of normally
distributed random vectors are still normally distributed.

In general, let $\mathbf{Y}=C\mathbf{X}$ be a $N$-dimensional random
vector with multi-dimensional distribution $N(\mu_Y,\Sigma_Y)$
written as a $N\times M$ linear transformation $C$ of a
$M$-dimensional random vector $\mathbf{X}$ with multi-dimensional
distribution $N(\mu_X,\Sigma_X)$ then:
\begin{equation}
    \Sigma_Y=C\Sigma_X C^T.\label{5.4}
\end{equation}
\noindent This result provides an easy way to generate a vector of
dependent normal random variables $\mathbf{Y}=C\mathbf{X} \sim
N(\mu_Y,\Sigma_Y)$ from a set of independent ones $X$. Indeed, the
dependence is completely taken into account by the covariance
matrix:
\begin{equation}
    \Sigma_Y = CC^T\label{5.5}
\end{equation}
\noindent The general problem consists of finding the linear
transformation $C$, (for further details and proofs see
Cufaro-Petroni \cite{Cufaro1996}).
\subsection{Cholesky Construction}
As far as the generation of a brownian motion is concerned, we note
that method $(\ref{5.1})$ can be written as:
\begin{equation}\label{5.1.1}
\left(
\begin{array}{c}
W_i\left( t_1\right) \\
\vdots \\
W_i\left( t_N\right)
\end{array}
\right) =C_R\left(
\begin{array}{c}
\epsilon_{1} \\
\vdots \\
\epsilon_N
\end{array}%
\right) \text{,}
\end{equation}
\noindent where $C_R$ is the Cholesky matrix associated to the
autocorrelation matrix of each brownian motion $W_i(t)$.

Referring to the general problem the Cholesky decomposition simply
faces the question of finding a matrix fulfilling equation
(\ref{5.5}) among all lower triangular matrices.

This is not a unique possibility, there are several other choices,
but all of them must satisfy the general problem (\ref{5.5}). We
will concentrate on two of them: the Principal Component Analysis
(PCA) proposed by Acworth, Broadie, and Glasserman (1998)
\cite{ABG1998} and a Kronecker Product Approximation that we
introduce as a different and new approach in Section 5.4.

We apply the Cholesky decomposition method in order to draw the
random vector $\mathbf{\epsilon}$ with distribution
$N(0,\Sigma_{MN})$.

In  case of constant volatilities we showed that $\Sigma_{MN} =
R\otimes \Sigma$. We can exploit the Kronecker product compatibility
with Cholesky decomposition to get:
\begin{equation}\label{5.1.2}
    C_{\Sigma_{MN}}=C_{R}\otimes C_{\Sigma}
\end{equation}
\noindent where $C_R$ is given by equation ($\ref{4.3}$). By means
of the Kronecker product we can reduce the computational effort by
splitting
 the analysis of an $MN\times MN$ matrix into the analysis
of two smaller $M\times M$ and $N\times N$ matrices .

When time-dependent volatilities are considered  we cannot exploit
the properties of the Kronecker product. $\Sigma_{MN}$ can be
partitioned into block matrices $\Sigma (t_1),\dots,\Sigma(t_N)$
that are not constant anymore and depend on the point of the time
grid.

Provided this time-dependent feature, all the information carried
out by $\Sigma_{MN}$ hinges in $N$ smaller $M\times M$ matrices.
These latter matrices depend on the particular time-dependent
functions that determine the evolution of the volatilities and on
the constant correlation among the assets returns (the analysis can
be applied to time-dependent instantaneous correlations).

In the following we present a faster than the standard Cholesky
decomposition algorithm that focuses on particular form of the
covariance matrix $\Sigma_{MN}$.

In the time-dependent volatility case the global covariance matrix
$\Sigma_{MN}$ satisfies the \emph{boomerang} shape property as $R$
as well as their Cholesky matrices. We consider this feature with
respect to the partitioned matrix notation.

It is possible to develop all the calculations storing $N$ block
matrices, $(\Sigma(t_1),\dots,\Sigma(t_N))$, in a tri-linear tensor
$(\Sigma_{tot})_{ikl}$. For any fixed $\hat{l}$ the block
$(\Sigma_{tot})_{ik\hat{l}}$ coincides with $\Sigma(t_{\hat{l}})$.
Consequently we perform the \emph{ad hoc} Cholesky decomposition
suited for partitioned \emph{boomerang} shaped matrices.

Using the partitioned matrix notation, the Cholesky algorithm develops according the following steps:
\begin{equation*}
\Sigma_{MN} = \left(\begin{array}{c|c} \Sigma_{TL} &
\Sigma_{BL}^T \\ \hline  \Sigma_{BL} &
\Sigma_{BR}\end{array}\right) = \left(\begin{array}{c|c}
C_{TL} & 0 \\ \hline C_{BL} &
 C_{BR}\end{array}\right)\left(\begin{array}{c|c}
 C^T_{TL} &  C^T_{BL}\\ \hline  0 &
C^T_{BR}\end{array}\right)
\end{equation*}
The block matrices with index \emph{TL} (Top-Left) are $M\times
M$, those ones with \emph{BL} (Bottom-Left) are $(N-1)M\times M$,
those ones with \emph{BR}(Bottom-Right) are $(N-1)M\times (N-1)M$.
\begin{enumerate}
\item{Decompose the Top-left block. \\
\begin{equation*}
\Sigma_{TL}=C_{TL}C_{TL}^T=C_1C_1^T
\end{equation*}}
\item{Decompose the Bottom-left block.
\begin{equation*}
\Sigma_{BL} = C_{BL}C_{TL}^T
\end{equation*}
In particular exploiting the \emph{boomerang} shape property we should have:\\
\begin{displaymath}
\Sigma_{BL} = \left(\begin{array}{c} \Sigma_{TL} \\ \vdots \\
\Sigma_{TL}\end{array}\right) = \left(\begin{array}{c}
C_{TL}C_{TL}^T
\\ \vdots \\ C_{TL}C_{TL}^T\end{array}\right)
\end{displaymath}
Due to the \emph{boomerang} shape structure of the global covariance
matrix, this second step can be avoided, because it consists of
repeating the first step.}
\item
The Cholesky decomposition is iterated to Bottom-Right
block.\\
\begin{equation*}
\Sigma_{BR} =
C_{BR}C_{BR}^T+C_{BL}C_{BL}^T
\end{equation*}
The last term on the right hand side of the previous equation is
known, because it has been calculated in step 1.

We let $\Sigma_{Update}$ define a $(N-1)M\times (N-1)M$ matrix by
the following expression:
\begin{displaymath}
\Sigma_{Update} = \Sigma_{BR} -
C_{BL}C_{BL}^T =
C_{BR}C_{BR}^T
\end{displaymath}
\noindent we can conclude that after decomposing $\Sigma_{Update}$
and getting $C_{BR}^T$ we have the complete picture of the global
Cholesky matrix.

This last step can be specified in greater detail referring to the \emph{boomerang}
shape feature of $\Sigma_{Update}$:
\begin{equation*}
\Sigma_{Update} = \left(\begin{array}{c|c} \Sigma(t_2) &
\tilde{\Sigma}_{TR} \\ \hline  \tilde{\Sigma}_{BL} &
\tilde{\Sigma}_{BR}\end{array}\right) - \left(\begin{array}{c}
C_{TL}\dots C_{TL}\end{array}\right)\left(\begin{array}{c}
C_{TL}^T \\ \vdots \\ C_{TL}^T\end{array}\right)
\end{equation*}
\noindent where $\tilde{\Sigma}_{BL}$ and $\tilde{\Sigma}_{BR}$, are
$(N-1)M \times M$ and $(N-1)M\times(N-1)M$ matrices. \noindent After
all the calculation we obtain:
\begin{equation*}
\Sigma_{Update}
=
\left(\begin{array}{c|c} \Sigma(t_2) - C_{TL}C_{TL}^T&
TR \\ \hline  BL &
BR\end{array}\right)
=C_{BR}C_{BR}^T
=
\left(\begin{array}{c|c}
C_{2} C_{2}^T & TR \\ \hline  BL &
BR\end{array}\right)
\end{equation*}
\noindent where $TR$, $BL$ and $BR$ are partitioned \emph{boomerang}
shaped matrices. $C_{2}$ represents the $M\times M$ Top-Left block
of $C_{BR}$, while $\Sigma(t_1) =C_{TL}C_{TL}^T=C_1C_1^T$
\end{enumerate}

\noindent The algorithm can be implemented running a loop of $N$ iterations.

The first iteration consists of realizing the Cholesky decomposition of step 1 described above.

The generic iteration $i$ consists in subtracting the Top-left Block
of the $i-1$ updated matrix to all the remaining $N-i$
 blocks (their dimension is  $M\times M$) of the tri-linear tensor $(\Sigma_{tot})_{ikl}$
 and that calculate the calculate the Cholesky decomposition.

This algorithm returns $N$ block matrices, whose dimension is
$M\times M$, that are stored in tri-linear  tensor,
$(C_{tot})_{ikj}$ that represents the global Cholesky
 matrix.

\subsection{Principal Component Analysis}
A more efficient approach for the path generation is based on
the Principal Component Analysis (PCA).

$\Sigma_{Y}$ is a symmetric matrix and can be diagonalized as
\begin{equation}\label{5.2.1}
   \Sigma_{Y} = E\Lambda E^T = (E \Lambda^{1/2}) (E \Lambda^{1/2})^T.
\end{equation}
\noindent For this method, the linear transformation $C$ solving
equation (\ref{5.5}) is defined as $E\Lambda^{1/2}$. $\Lambda$ is
the diagonal matrix of all the positive eigenvalues of $\Sigma_{Y}$
sorted in decreasing order  and $E$ is the orthogonal matrix
($EE^T=I$) of all the correspondent eigenvectors.

The matrix $E\Lambda^{1/2}$ has no particular structure and
generally does not provide computational advantage with respect to
the Cholesky decomposition.

This transformation can be interpreted as a sort of rotation of
the random vector whose covariance matrix is $\Sigma_{Y}$; in the
new frame of reference it has independent components whose
variances are the elements on the diagonal of $\Lambda$.

The higher efficiency of this method is due to the statistical
interpretation of the eigenvalues and eigenvectors (see Glasserman
\cite{Glass2004}).

Suppose we want to generate $\mathbf{Y}\sim N(0,\Sigma_Y)$ from a vector
$\mathbf{\epsilon} \sim N(0,I)$, we know that the random vector can be set as:
\begin{equation*}
    \mathbf{Y} = \sum_{k=1}^d \mathbf{c_k} \epsilon_k
\end{equation*}
\noindent where $\mathbf{c_k}$ is the $k$-th column of $C$.

Assume $\Sigma_Y$ has full rank $d$, then it is non singular and
invertible and the factors $\epsilon_k$ are themselves linear
combination of $Y_k$. In the special case $C = E\Lambda^{1/2}$,
$\epsilon_k$ is proportional to $\mathbf{e_k}\cdot \mathbf{Y}$.

The factors $\epsilon_k$ constructed in the previous way are optimal in a precise statistical sense.

Suppose we want to find the best singled-factor approximation of
$\mathbf{Y}$, that is to find the best linear approximation that
best captures the variability of the components of $\mathbf{Y}$. The
optimization problem consists in maximizing the variance of
$\mathbf{w}\cdot \mathbf{Y}$ with constraint of the form
$\mathbf{w}\cdot \mathbf{w}=1$:
\begin{equation}\label{5.2.2}
  \max_{\mathbf{w} \cdot \mathbf{w} =1} \quad \mathbf{w}\cdot \Sigma_Y \mathbf{w}
\end{equation}
\noindent If we sort the eigenvalues of $\Sigma_Y$ in decreasing
order then the optimization problem is solved by $\mathbf{e_1}$.
More generally  the best $k$-factors approximation of $\mathbf{Y}$
leads to factors proportional to $\mathbf{e_1}\cdot
\mathbf{Y},\dots,\mathbf{e_k}\cdot \mathbf{Y}$ with
$\mathbf{e_l}\cdot \mathbf{e_m} = \delta_{lm}$, with:
\begin{equation}
    \epsilon_k = \frac{1}{\sqrt{\lambda_k}}\mathbf{e_k} \cdot \mathbf{Y}.\label{5.2.3}
\end{equation}
This representation can be recasted as the minimization of the
mean squared error:
\begin{equation}\label{5.2.4}
    MSE = \mathbb{E}\left[ \|\mathbf{Y} - \sum_{i=1}^k \mathbf{c_i} \epsilon_i \|^2\right]
\end{equation}
\noindent where we are looking for the best $k$-factors mean
square approximation of $X$. This formulation gives the same
results.

In the statistic literature the linear combination
$\mathbf{e_k}\cdot \mathbf{Y}$ is called principal component of
$\mathbf{Y}$. The amount of variance explained by the first $k$
principal components is the ratio:
\begin{equation}
    \frac{\sum_{i=1}^k \lambda_i}{\sum_{i=1}^d \lambda_i}\label{5.2.5}
\end{equation}
\noindent where $d$ is the rank of $\Sigma_Y$.

We can apply PCA to generate a one-dimensional brownian motion BM$(0,R)$
calculating the eigenvectors and eigenvalues of the sampled auto-covariance matrix $R$ and then rearranging them in decreasing order.
The magnitude of the eigenvalues of this matrix drops off rapidly.
For instance it is possible to verify that in the case of a brownian motion with $32$ time steps the amount of variance explained by the first
five factors is $81\%$ while it exceeds $99\%$ at $k=16$.

This result is fundamental in identifying the effective dimension of
the integration problem. PCA helps Monte Carlo estimation procedures
based on the generation of brownian motion where we should identify
the effective dimension of the problem. With this choice we can
identify the most important factors in a precise statistical
framework by fixing a value $p$ in the determining the effective
dimension. (for instance $p=99\%$).

This statistical ranking of the normal factors cannot be implemented
by Cholesky decomposition that we expect will return unbiased Monte
Carlo estimations but higher RMSEs.

As far as the multi-dimensional brownian motion is concerned, we start with the constant volatility case.
We have already shown in section 4 that the covariance matrix $\Sigma_{MN}$ of the multi-dimensional
brownian motion BM$(0,\Sigma_{MN})$ can be written as $R\otimes \Sigma$.

Property 10 of the Kronecker product permits to improve the speed of
the computation of the eigenvalues and eigenvectors of
$\Sigma_{MN}$. It reduces this calculation into the computation of
the eigenvalues and vectors of the two smaller matrices $R$ and
$\Sigma$.

Coupling the use of the Kronecker product analysis with the ANOVA
definition of effective dimension we can implement a fast and
efficient Monte Carlo estimation in order to price exotic
multi-dimensional path-dependent options.

Empirical evidence in finance shows that effective dimension is
often lower than the problem dimension $d$, (see Caflisch,
Morokoff, and Owen \cite{CMO1997} for a general discussion). We
focus our analysis on Asian options pricing after formulating the
pricing problem an integral. As we presented in section 3 ANOVA is
used as to provide a representation of the integrand
 as a sum of orthogonal functions. If each of these orthogonal functions depends only on a distinct subset of the coordinates,
 the integrand can be written as a sum of integrals of functions of lower dimension.
  The complexity of the computation of the integral has been reduced with respect to the integral dimension.
In pricing Asian options we are not able to reduce the dimension of
the original integrand by this approach, because we cannot exactly
find a set of orthogonal functions. What we can propose is an
approximation based on the PCA construction. In our finance problems
we achieved a representation involving matrices,  describing the
dependence between the different variables, as arguments of the
exponential function $g(\centerdot)$. Our approximation consists in
a direct application of ANOVA and effective dimension calculation to
the random vector $\mathbf{Z}$. This is equivalent to the Taylor
expansion up to the first order of the exponential function
$g(\centerdot)$ that leads to the following definition of effective
dimension, $d_T$, of the problem (in truncation sense):
\begin{equation}
\sum_ {d=1}^{d_T} \lambda_d \leq tr(\Lambda)p\label{5.2.6}
\end{equation}
where $\lambda_d \in \sigma(\Sigma_{MN})$. The level $p$ is
arbitrary; we chose $p=99\%$.
\subsection{The Kronecker Product Approximation}
The time-dependent volatilities market has a covariance matrix
$\Sigma_{MN}$ with time-dependent blocks. Generally, it has not a
particular expression because it depends on the volatility functions
and the instantaneous correlation. The covariance matrix of the
asset returns is not anymore constant so that $\Sigma_{MN}$ cannot
be written as a Kronecker product.

We have shown that a fast Cholesky decomposition algorithm can be
ran but it does not take any ANOVA and effective dimension
consideration, while the PCA approach is still applicable but we
cannot reduce the computational burden using the properties offered
by the Kronecker product.

In the constant volatility case the special structure of
$\Sigma_{MN}$ makes  possible to compute all the eigenvalues and
eigenvectors with $M^3 + N^3$ operations, written $O(M^3 + N^3)$,
instead of $O\left((MN)^3\right)$ for a general $MN\times MN$ square
matrix.

The market under consideration has the multi-dimensional brownian
motion as unique source of risk. Its generation procedure is
independent of the constant or time-dependent volatilities because
its autocovariance matrix $R$ is not influenced by these market
features.

Based on these considerations our proposition is to find a constant
covariance matrix among the assets, $K$, in order to approximate, in
an appropriate sense, the global covariance matrix $\Sigma_{MN}$ as
a Kronecker product of $R$ and $K$. Our hypothesis is that the
effective dimension of the problem should not dramatically change
after this transformation with an advantage from the computational
point of view. We develop the PCA decomposition of the approximating
matrix assuming that the principal components are not so different
from those of the original random vector. This approximation would
lead to a different multi-dimensional path because $R\otimes K$ is
not the covariance matrix of the original process. The global and
true path is reobtained using the Cholesky factorization.

In the following we illustrate the proposed procedure that we label KPA.

The general problem consists of finding two matrices
$B\in\mathbb{R}^{m_1\times n_1}$ and $C\in\mathbb{R}^{m_2\times
n_2}$ that minimize the Frobenius norm. All calculations and
proofs can be found in Pitsianis, Van Loan \cite{PV1993} and Van
Loan \cite{VL2000}:
\begin{equation}
\Phi_A(B,C) = \parallel A - B\otimes C \parallel^2\label{5.4.1}
\end{equation}
where $A\in\mathbb{R}^{m\times n}$ is an assigned matrix with
$m=m_1m_2$ and $n=n_1n_2$.
\\ The main idea  is to look for a rearrange
matrix $\mathcal{R}(A)$ such that equation (\ref{5.4.1}) can be
rewritten as $\Phi_A(B,C) = \parallel \mathcal{R}(A) - vec(B)\otimes
vec(C)^T
\parallel^2$.
\newtheorem{vecu}[Definition1]{Definition}
\begin{vecu}[\textbf{The vec operation}]
 The vec operation transforms a matrix
$X\in\mathbb{R}^{M,N}$ into a column vector
$vec(X)\in\mathbb{R}^{MN}$ by 'stacking' the columns:
\end{vecu}
\begin{displaymath}
A=\left(\begin{array}{cc}a_{11} & a_{12}\\a_{21} &
a_{22}\end{array}\right)\Longrightarrow vec(X) =
\left(\begin{array}{c}a_{11} \\ a_{21}\\a_{12}
\\a_{22}\end{array}\right)
\end{displaymath}

As far as our approximation is concerned the general problem is
simplified. Indeed, the new problem consists of finding only one
matrix $K$ minimizing the Frobenius norm:
\begin{equation}\label{5.4.2}
\Phi(K) = \parallel \Sigma_{MN} - R\otimes K
\parallel^2
\end{equation}
\noindent The approach is equivalent to a Least Square problem in the
${K_{ik}}$.

The elements $K_{ik}$ are given by the formula below (for a complete
proof see Pitsianis, Van Loan \cite{PV1993} p. 8):
\begin{equation}
{K_{ik}} = \frac{tr\Big(\mathcal{R}(
\Sigma_{MN})_{ik}R\Big)}{tr\Big(R
R^T\Big)}\label{5.4.3}
\end{equation}
where $\mathcal{R}({\mathbf \Sigma_{MN}}))_{ik}$ is a $N\times N$
matrix. For any $i$ and $k$ ranging from $1$ to $M$,
$\mathcal{R}({\mathbf \Sigma_{MN}}))_{ik}$ is obtained by sampling
$\Sigma_{MN}$ with $M$ as sampling step.

By its definition,  it can be noticed that for any $i$ and $k$
$\mathcal{R}(\Sigma_{MN}))_{ik}$ is a  \emph{boomerang} shaped block
matrix.

By direct computations and relying on the particular form of
$R$, the denominator of the equation(\ref{5.4.3}) is:
\begin{equation}
tr\Big(RR^T\Big) = tr\Big(R^2\Big) = \sum_{j=1}^N
\Big(2(N-j)+1\Big)t_j^2\label{5.4.4}
\end{equation}

Moreover, given two general $N\times N$ \emph{boomerang} shaped
matrices $A$ and $B$ the
trace of their product is:
\begin{equation}
tr(A^TB) = tr(AB) = \sum_{j=1}^N
\Big(2(N-j)+1\Big)a_{jj}b_{jj}\label{5.4.5}
\end{equation}
$a_{jj}$ and $b_{jj}$ are the only significant value to store.

The considerations above permit to evaluate $K$ in a fast and
efficient way without high computational efforts.

As already mentioned, if we would use the ANOVA-PCA procedure to
$R\otimes{K }$ we would not get the required path. Let $E$ and
$\Lambda$ be the eigenvectors and eigenvalue matrices associated to
$R\otimes{K }$, if we would consider $ E\Lambda^{1/2}$ as a
generating matrix we would generate a path whose global covariance
matrix is $R\otimes{K }$ and not $ \Sigma_{MN}$.

In order to tackle to the original problem the Cholesky
decomposition is used. In fact given two $N$ dimensional random
vectors $\mathbf{Z_1}$ and $\mathbf{Z_2}$ with covariance matrices $\Sigma_1$ and
$\Sigma_2$ respectively, we can always write:
\begin{equation}\label{5.4.6}
\bigg\{
\begin{array}{c}
  \mathbf{Z_1} = C_1 \mathbf{\epsilon} \\
  \mathbf{Z_2} = C_2 \mathbf{\epsilon}
\end{array}
\end{equation}
\noindent where $C_1$ and $C_2$ are the Cholesky matrices of
$\Sigma_1$ and $\Sigma_2$, respectively and $\mathbf{\epsilon}$ is a
vector of independent random variables. At the same time we can
generate $Z_2$ by PCA:
\begin{equation}\label{5.4.7}
    \mathbf{Z_2} = E_2\Lambda_2^{1/2}\mathbf{\epsilon}
\end{equation}
\noindent where $E_2$ and $\Lambda_2$ comes
from the complete PCA of $\Sigma_2$.

Combining the above equalities we have:
\begin{equation}
\mathbf{Z_1} = C_1C_2^{-1}E_2\Lambda_2^{1/2}\mathbf{\epsilon}\label{5.4.8}
\end{equation}
It is possible to generate a random path $\mathbf{Z_1}$ applying the
$PCA$ to $\mathbf{Z_2}$ and than turning back to the original
problem. Our fundamental assumption is that the effective dimension
of our problem remains almost unchanged and, in the estimation
procedure, we apply almost the same statistical importance to the
original principal components giving an advantage from the
computational point of view.

Focusing this result to the problem under study, we let
$\Sigma_1=\Sigma_{MN}$ and $\Sigma_2 =R\otimes K $ so that
equation(\ref{5.4.8}) becomes:
\begin{equation}\label{5.4.9}
\mathbf{Z} = C_{\Sigma_{MN}}C_R^{-1}\otimes C_K^{-1}
E_2\Lambda^{1/2}\mathbf{\epsilon}
\end{equation}
$C_{\Sigma_{MN}}$, $ C_R^{-1}$ and $C_K^{-1}$ are the Cholesky
matrices of $\Sigma_{MN}$, $R$ and $K$, respectively. In the
derivation of the previous equation we exploit several properties of
the Kronecker product.

We again stress the fact that in the case of time-dependent
volatilities we analyze the effective dimensions of the integral
problem after the Kronecker product approximation. Generally this
second approximation would return a higher effective dimension with
respect to the normal case where only a linear approximation is
considered. Furthermore our method generates the correct required
path as proved.

In order to obtain a fast and efficient algorithm for the path
generation we develop all the calculations:
\begin{enumerate}
\item{$C_R^{-1}\otimes C_K^{-1}$ is a sparse
bi-diagonal partitioned matrix:
\begin{displaymath}
C_R^{-1}\otimes C_K^{-1} =
\left(\begin{array}{ccccc} \frac{C_K^{-1}}{\sqrt{t_{1}}}
& 0 & \ldots & \ldots & 0
\\ -\frac{C_K^{-1}}{\sqrt{t_2-t_1}} & \frac{C_K^{-1}}{\sqrt{t_2-t_1}} & 0 & \dots & 0
\\ 0 & -\frac{C_K^{-1}}{\sqrt{t_3-t_2}} & \frac{C_K^{-1}}{\sqrt{t_3-t_2}} &
\ddots & \vdots
\\ \vdots & \vdots & \ddots & \ddots & \vdots
\\ 0 & 0 & 0 & -\frac{C_K^{-1}}{\sqrt{t_n-t_{n-1}}} &
 \frac{C_K^{-1}}{\sqrt{t_n-t_{n-1}}} \end{array}\right)
\end{displaymath}}

\item {$C_{\Sigma_{MN}}C_R^{-1}\otimes C_K^{-1}$ is lower triangular partitioned
matrix.\\ \\ $C_{\Sigma_{MN}}C_R^{-1}\otimes
C_K^{-1}=$
\begin{displaymath}
\left(\begin{array}{ccccc}
\frac{C_{\Sigma1}}{\sqrt{\Delta_{1}}}C_K^{-1} & 0 &
\ldots & \ldots & 0
\\ \Big(\frac{C_{\Sigma 1}}{\sqrt{\Delta_1}}-\frac{\mathbf{C_{\Sigma
2}}}{\sqrt{\Delta_2}}\Big)C_K^{-1} &
\frac{C_{\Sigma 2}}{\sqrt{\Delta_2}}C_K^{-1} &
0 & \dots & 0
\\ \vdots & \Big(\frac{C_{\Sigma 2}}{\sqrt{\Delta_2}}-\frac{C_{\Sigma
3}}{\sqrt{\Delta_3}}\Big)C_K^{-1} &\qquad
\frac{C_{\Sigma 3}}{\sqrt{\Delta_3}}C_K^{-1} &
\ddots & \vdots
\\ \vdots & \vdots & \ddots & \ddots & \vdots
\\ \Big(\frac{C_{\Sigma 1}}{\sqrt{\Delta_1}}-\frac{C_{\Sigma
2}}{\sqrt{\Delta_2}}\Big)C_K^{-1} &
\Big(\frac{C_{\Sigma
2}}{\sqrt{\Delta_2}}-\frac{C_{\Sigma
3}}{\sqrt{\Delta_3}}\Big)C_K^{-1} & \dots &
\Big(\frac{C_{\Sigma
N-1}}{\sqrt{\Delta_{N-1}}}-\frac{C_{\Sigma
N}}{\sqrt{\Delta_N}}\Big)C_K^{-1} &
 \frac{C_{\Sigma N}}{\sqrt{\Delta_N}}{C_K^{-1}} \end{array}\right)
\end{displaymath}}
\end{enumerate}
$C_{\Sigma i}$ for $i=1,\dots,N$ indicates the $i$-th
block matrix of the tri-linear tensor $(C_{tot})_{ikj}$.
$\Delta_{i}=t_i-t_{i-1}$ where $t_0=0$ is understood.

Only $( C_{tot})_{ikj}$ and the sequence
$\{\Delta_i\}_{i=1,\dots,N}$ need to store all the information
embedded in $C_{\Sigma_{MN}}C_R^{-1}\otimes C_K^{-1}$.\\ The total
generating matrix $C_{\Sigma_{MN}}C_R^{-1}\otimes
 C_K^{-1} E_2\Lambda^{1/2}$ can be computed quickly by
matrix product with partitioned matrices.
\section{Solution Methodology}

\noindent We aim to provide an efficient technique that improves the
precision of the general Monte Carlo method to exotic derivative
contracts and in particular Asian options. According to equation
(\ref{2.2.10}) the actual problem consists of generating a sample of
uniform random draws to uniformly cover the whole hypercube $\left[
0,1\right] ^{d}$. In the following subsections we introduce
different ways to generate random numbers that uniformly cover the
hypercube $[0,1]^d$.
%
%
\subsection{Stratification and Latin Hypercube Sampling}

\noindent Stratified sampling is a variance reduction method for Monte Carlo
estimates. It amounts to partitioning the hypercube $\mathcal{D}=\left[
0,1\right) ^{d}$ into $H$ disjoint strata $\mathcal{D}_{h}$, ($h=1,\dots ,H$%
), \textit{i.e.}, $\mathcal{D}=\bigcup_{i=1}^{H}\mathcal{D}_{h}$ where $%
\mathcal{D}_{k}\bigcap \mathcal{D}_{j}=\varnothing $ for all $j\neq
k$, then estimating the integral over each set, and finally summing
up these numbers (see Boyle, Broadie and Glasserman \cite{BBG1997}
for more on this issue). Specifically, mutually independent uniform
samples $x_{1}^{h},\dots ,x_{n_{h}}^{h}$ are simulated within a
stratum $\mathcal{D}_{h}$, and the resulting integrals are combined.
The resulting stratified sampling estimator is unbiased. Indeed:
\begin{eqnarray*}
\mathbb{E}\left[ \widehat{I}_{strat}\right] &=&\sum_{h=1}^{H}\frac{|\mathcal{%
D}_{h}|}{n_{h}}\sum_{i=1}^{n_{h}}\mathbb{E}\left[ f\left( x_{i}^{h}\right) %
\right] \\
&=&\sum_{h=1}^{H}|\mathcal{D}_{h}|\mu _{h} \\
&=&\sum_{h=1}^{H}\int_{\mathcal{D}_{h}}f\left( x\right) \,dx=I.
\end{eqnarray*}%
\noindent where $|\mathcal{D}_{h}|$ denotes the volume of stratum $D_{h}$.
Moreover, this estimator displays a lower variance compared to a crude Monte
Carlo estimation, \textit{i.e.},%
\begin{equation*}
\text{Var}\left[ \widehat{I}_{strat}\right] \leq \frac{\sigma ^{2}}{n}.
\end{equation*}%
\noindent Stratified sampling transforms each uniformly distributed
sequence $\mathbf{U}_{j}=\left( U_{1j},\dots ,U_{dj}\right) $ in
$\mathcal{D}$ into a new
sequence $\mathbf{V}_{j}=\left( V_{1j},\dots ,V_{dj}\right) $ according to the rule%
\begin{equation*}
\mathbf{V}_{j}=\frac{\mathbf{U}_{j}+\left( i_{1},\dots ,i_{d}\right)
}{n},j=1,\dots ,n,i_{k}=0,\dots ,n-1,k=1,\dots ,d.
\end{equation*}%
where $\left( i_{1},\dots ,i_{d}\right) $ is a deterministic
permutation of the integers $1$ through $d$. This procedure ensures
that one $\mathbf{V}_{j}$ lies in each of the $n^{d}$ hypercubes
defined by the stratification.
\begin{figure}
\begin{center}
\includegraphics[width=0.7\textwidth]{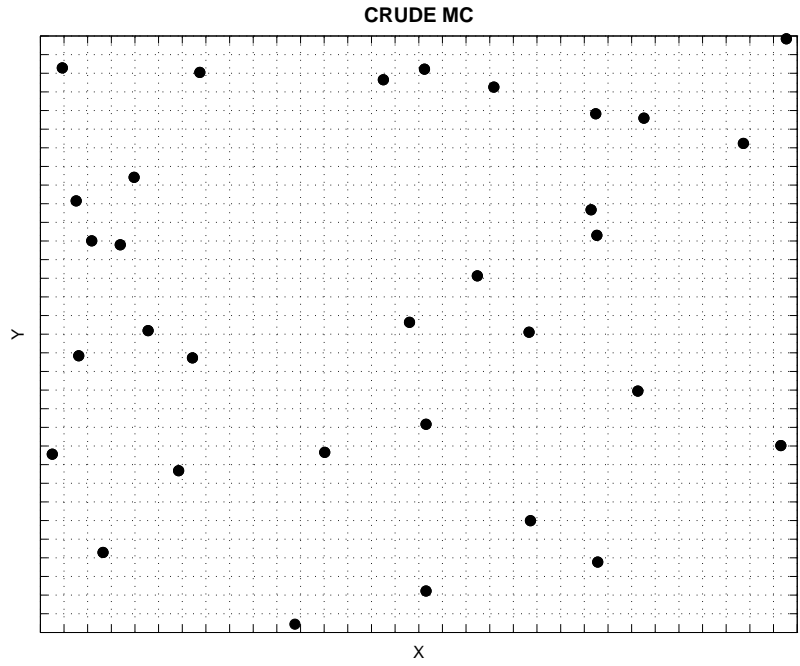}
\end{center}
\caption{The panel shows 32 points drawn with standard pseudorandom
generators}
\end{figure}
Latin Hypercube Sampling (LHS) can be seen as a way of randomly
sampling $n$ points of a stratified sampling while preserving the
regularity from stratification (see, for instance, Glasserman
\cite{Glass2004}). Let $\pi _{1},\dots ,\pi _{d}$ be independent
random permutations of the first $n$ positive integers, each of
them uniformly distributed over the $n!$ possible permutations.
Set
\begin{equation}\label{6.1.1}
T_{jk}=\frac{U_{jk}+\pi _{k}\left( j\right) -1}{n},\qquad j=1,\dots
,n,k=1,\dots ,d,
\end{equation}%
\noindent

\noindent where $\pi _{k}\left( j\right) $ represents the $j$-th
component of the permutation for the $k$-th coordinate.
Randomization ensures that each vector $\mathbf{T}_{j}$ is uniformly
distributed over the $d$ dimensional hypercube. Moreover, all
coordinates are perfectly stratified since there is exactly one
sample point in each hypercube of volume $1/n$. For $d=2$, there is
only one point in the horizontal or vertical stripes of surface
$1/n$
(see Figure 2). The base and the height are $1/n$ and 1, respectively. For $%
d>2$ it works in the same way. It can be proven that for all $n\geq
2,d\geq 1 $ and squared integrable functions $f$, the error for the
estimation with the Latin Hypercube Sampling is smaller or equal to
the error for the crude Monte Carlo (see Koehler and Owen
\cite{KO1996}):
\begin{equation}
Var\left[ \widehat{I}_{LHS}\right] \leq \frac{\sigma ^{2}}{n-1}.\label{6.1.2}
\end{equation}%
Figure 2 shows the distribution of 32 points generated with the LHS
method. For the LHS method we notice that there is only $1$ point
(dotted points in Figure 2) in each vertical or horizontal stripe
whose base is $1$ and height is $1/32$: it means that there is only
a vertical and horizontal stratification.
\begin{figure}\label{fig1}
\begin{center}
\includegraphics[width=0.7\textwidth]{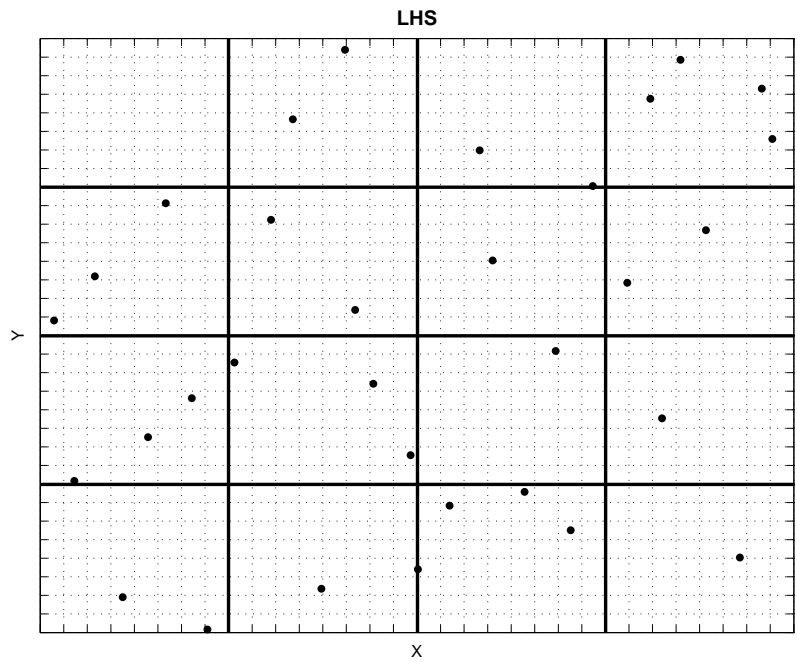}
\end{center}
\caption{The panel shows 32 points generated with LHS}
\end{figure}

\subsection{Low-Discrepancy Sequences}

\noindent As previously mentioned, the standard \emph{MC} method is
based on a completely random sampling of the hypercube $\left[
0,1\right) ^{d}$ and its precision can be improved using
stratification or Latin Hypercube sampling. These two methods ensure
that there is only one point in each smaller hypercube fixed by the
stratification as illustrated in Figure 2.
 At the same time, these techniques provide nothing more
than the generation of uniform random variables in smaller sets.

A completely different way to approach the sampling problem is to
build-up a deterministic sequence of points that uniformly covers
the hypercube $\left[ 0,1\right) ^{d}$ and to run the estimation
using this sequence. Obviously, there is no statistical quantity
that may represent the uncertainty since the estimation always gives
the same results. The Monte Carlo method implemented with the use of
low-discrepancy sequences is called Quasi-Monte Carlo (QMC).

The mathematics involved in generating a low-discrepancy sequence is
complex and requires the knowledge of the number theory. In the
following, only an overview of the fundamental results and
properties is presented (see Niederreiter \cite{Ni1992} for more on
this issue).

We define the quantity $D_{n}^{\ast }=D_{n}^{\ast }\left( P_{1},\dots
,P_{n}\right) $ as the star discrepancy. It is a measure of the uniformity
of the sequence $\left\{ P_{n}\right\} _{n\in \mathbb{N}^{\ast }}\in \left[
0,1\right) ^{d}$ and it must be stressed that it is an analytical quantity
and not a statistical one. For example, if we consider the uniform
distribution in the hypercube $\left[ 0,1\right) ^{d},$ the probability of
being in a subset of the hypercube is given by the volume of the subset. The
discrepancy measures how the pseudo-random sequence is far from the
idealized uniform case, i.e. it is a measure, with respect to the $L_{2}$
norm for instance, of the inhomogeneity of the pseudo-random sequence.
\newtheorem{Low}[Definition1]{Definition}
\begin{Low}[\textbf{Low-Discrepancy Sequencies}]
A sequence $\left\{ P_{n}\right\} _{n\in \mathbb{N}^{\ast }}$ is
called low-discrepancy sequence if:
\end{Low}
\begin{equation}
D_{n}^{\ast }\left( P_{1},\dots ,P_{n}\right) =O\left( \frac{\left( \ln
n\right) ^{d}}{n}\right) .  \label{6.2.1}
\end{equation}

\noindent i.e. if its\textbf{\ }star discrepancy decreases as $\left( \ln
n\right) ^{d}/n$\textbf{.}

The following inequality, attributed to Koksma and Hlawka, provides
an upper bound to the estimation error of the unknown integral with
the QMC method in terms of the star discrepancy:
\begin{equation}
|I-\hat{I}|\leq D_{n}^{\ast }\,\,V_{HK}\left( f\right) .  \label{6.2.2}
\end{equation}%
\noindent $V_{HK}\left( f\right) $ is the variation in the sense of
Hardy and Krause. Consequently, if $f$ has a finite variation and
$n$ is large enough, the QMC approach gives an error smaller than
the error obtained by the crude \emph{MC} method for low dimensions
$d$. However, the problem is difficult owing to the complexity if
estimating the Hardy-Krause variation, which depends on the
particular integrand function.

In the following sections we briefly present digital nets and the
well-known Sobol' sequence that is the most frequently used low
discrepancy sequence to run Quasi-Monte Carlo simulations in
finance.
%
%
\subsection{Digital Nets}

\noindent Digital nets or sequences are obtained by the number
theory and owe their name to the fact that their properties can be
recognized by their digital $b$-ary expansion in base $b$. Many
digital nets exist; the ones most often used and considered most
efficient are the Sobol' and the Niederreiter-Xing sequences.

The first and simplest digital sequence with $d=1$ is due to Van der Corput
and is called the radical inverse sequence. Given an integer $b\geq 2$, any
non-negative number $n$ can be written in base $b\ $as:
\begin{equation}
n=\sum_{k=1}^{\infty }n_{k}b^{k-1}.\label{6.3.1}
\end{equation}%
\noindent The base $b$ radical inverse function $\phi _{b}\left( n\right) $
is defined as:
\begin{equation}\label{6.3.2}
\phi _{b}\left( n\right) =\sum_{k=1}^{\infty }n_{k}b^{-k}\quad \in \left[
0,1\right) ,
\end{equation}%
\noindent where $n_{k}\in \left\{ 0,1,\dots ,b-1\right\} $ (Galois set).

By varying $n$ the Van der Corput sequence is constructed.
\begin{table}[tbp] \centering%
\begin{tabular}{||l|l|l|l||}
\hline\hline
$N$& $n$ base 2 & $\phi _{2}\left( n\right) $
base 2 & $\phi _{2}\left( n\right) $ \\ \hline
0 & 000. & 0.000 & 0.000 \\
1 & 001. & 0.100 & 0.500 \\
2 & 010. & 0.010 & 0.250 \\
3 & 011. & 0.110 & 0.750 \\
4 & 100. & 0.001 & 0.125 \\
5 & 101. & 0.101 & 0.625 \\
6 & 110. & 0.011 & 0.375 \\
7 & 111. & 0.111 & 0.875 \\ \hline\hline
\end{tabular}%
\caption{Van der Corput sequence.}\label{VanCorp_S}%
\end{table}%
Table \ref{VanCorp_S} illustrates the first seven Van der Corput points for $%
b=2$. Consecutive integers alternate between odd and even; these points
alternate between values in $\left[ 0,1/2\right) $ and $\left[ 1/2,1\right) $%
. The peculiarity of this net is that any consecutive $b^{m}$ points from
the radical inverse sequence in base $b$ are stratified with respect to $%
b^{m}$ congruent intervals of length $b^{-m}$. This means that in each
interval of length $b^{-m}$ there is only one point.

Table \ref{VanCorp_S} shows an important property that is exploited
in order to generate digital nets, because a computing machine can
represent each number with a given precision, referred to as
\textquotedblleft machine epsilon\textquotedblright . Let
$z=0.z_{1}z_{2}\dots (base$ $b)\in \left[ 0,1\right) $ , define
$\mathbf{\Psi} (z)=(z_{1},z_{2},\dots )$ the vector of the its
digits, and truncate its digital expansion at the maximum allowed
digit $w$: $z=\sum_{k=1}^{\emph{w}}z_{k}b^{-k}$. Let $n=\left[
b^{w}z\right] =\sum_{h=1}^{w}n_{h}b^{h-1}\in N^{\ast }$, where
$\left[ x\right] $ denotes the greatest integer less than or equal
to $x$. It can be easily proven that:

\begin{equation*}
n_{h}=z_{w-h+1}\left( z\right) \quad \forall h=1,\dots ,w.
\end{equation*}%
\noindent This means that the finite sequences $\left\{
n_{h}\right\} _{h\in \left\{ 1,\dots ,w\right\} }$ and $\left\{
z_{k}\right\} _{k\in \left\{ 1,\dots ,w\right\} }$\ have the same
elements in opposite order. For example, in the table
\ref{VanCorp_S} we allow only $3$\ digits; in order to find the
digits of $\phi _{2}\left( 1\right) =0,5$ we consider $\phi
_{2}\left( 1\right) 2^{3}=4=0n_{1}+n_{2}0+n_{3}1$. The digits of
$\phi
_{2}\left( 1\right) $ are then $(1,0,0)$ as shown in the table \ref%
{VanCorp_S}.

The peculiarity of the Van der Corput sequence is largely required
in high dimensions, where the contiguous intervals are replaced by
multi-dimensional sets called  b-adic boxes.
\newtheorem{Boxi}[Definition1]{Definition}
\begin{Boxi}[\textbf{b-iadic Box}]
Let $b\geq 2$, $k_{j}$, $l_{j}$ with $0\leq l_{j}\leq b^{k_{j}}$ be all
integer numbers. The following set is called b-iadic box:
\end{Boxi}
\begin{equation}\label{6.3.3}
\prod_{j=1}^{d}\left[ \frac{l_{j}}{b^{k_{j}}},\frac{l_{j}+1}{b^{k_{j}}}%
\right) ,
\end{equation}%
\noindent where the product represents the Cartesian product.
\newtheorem{Net}[Definition1]{Definition}
\begin{Net}[\textbf{(t,m,d) Nets}]
Let $t\leq m$ be a non-negative integer. A finite set of points from $\left[
0,1\right) ^{d}$ is called $\left( t,m,d\right) $-net if every b-adic box of
volume $b^{-m+t}$ (bigger than $b^{-m}$) contains exactly $b^{t}$ points.
\end{Net}
This means that cells that \textquotedblleft should have\textquotedblright\ $%
b^{t}$ points do have $b^{t}$ points. However, considering the smaller
portion of volume $b^{-m}$, it is not guaranteed that there is just one
point.

A famous result of the theory of digital nets is that the
integration over a $\left( t,m,d\right) $ net can attain an
accuracy of the order of $O\left( \ln ^{d-1}\left( n\right)
/n\right) $ while, restricting to $\left( t,d\right) $ sequences,
it raises slightly to $O\left( \ln ^{d}\left( n\right) /n\right) $
(see Niederreiter \cite{Ni1992}). The above results are true only
for functions with bounded variation in the sense of Hardy-Krause.
%
%
\subsection{The Sobol' Sequence}
The Sobol' sequence is the first $d$ dimensional digital sequence, ($b=2$),
ever realized. Its definition is complex and is covered only briefly in the
following.
\newtheorem{Sob}[Definition1]{Definition}
\begin{Sob}[\textbf{The Sobol' Sequence}]
Let $\left\{ n_{k}\right\} _{k\in \mathbb{N^{\ast }}}$ be the
set of the b-ary expansion in base $b=2$ of any integer $n$; the $n$-th element $%
S_{n}$ of the Sobol' sequence is defined as:
\end{Sob}
\begin{equation}\label{6.4.1}
S_{n}= \sum_{k=1}^{+\infty }\left(n_{k}\,V_{k}\,mod\,
2\right)\,2^{-k},
\end{equation}%
\begin{table}[tbp] \centering%
\begin{tabular}{||l|l|l|l|l||}
\hline\hline d & P & M & Principal polynomial & q \\ \hline
1 & [1] & [1] & $1$ & 0 \\
2 & [1 1] & [1] & $x+1$ & 1 \\
3 & [1 1 1] & [1 1] & $x^{2}+x+1$ & 2 \\
4 & [1 0 1 1] & [1 3 7] & $x^{3}+x+1$ & 3 \\
5 & [1 1 0 1] & [1 1 5] & $x^{3}+x^{2}+1$ & 3 \\
6 & [1 0 0 1 1] & [1 3 1 1] & $x^{4}+x+1$ & 4 \\
7 & [1 1 0 0 1 ] & [1 1 3 7] & $x^{4}+x^{3}+1$ & 4 \\
8 & [1 0 0 1 0 1] & [1 3 3 9 9] & $x^{5}+x^{2}+1$ & 5 \\
9 & [1 1 1 0 1 1] & [1 3 7 13 3] & $x^{5}+x^{4}+x^{3}+x+1$ & 5 \\
10 & [1 0 1 1 1 1] & [1 1 5 11 27] & $x^{5}+x^{3}+x^{2}+x+1$ & 5 \\
\hline\hline
\end{tabular}%
\caption{Initial values satisfying Sobol´ property A up to dimension 10.
 By convention, the recurrence relation for the $0$-degree polynomial is $M_k\equiv1$}\label{DirNumb_S}%
\end{table}%
\noindent where $V_{k}\in \left[ 0,1\right) ^{d}$ are called
direction numbers. In practice, the maximum number of digits, $w$,
must be
given. In Sobol's original method the $i$-th number of the sequence $S_{ij}$%
, $i\in {\mathbb{N}},j\in \left\{ 1,\dots ,d\right\} $, is generated by
XORing (bitwise exclusive OR) together the set of $V_{kj}$ satisfying the
criterion on $k$ : the $k$-th bit of $i$ is nonzero. Antonov and Saleev
derived a faster algorithm by using the Grey code. Dropping the index $j$
for simplicity, the new method allows us to compute the $\left( i+1\right) $%
-th Sobol' number from the $i$-th by XORing it with a single
$V_{k}$, namely with $k$, the position of the rightmost zero bit
in $i$ (see, for instance, Press, Teukolsky, Vetterling and
Flannery \cite{PTVF1992}). Each different Sobol' sequence is based
on a different primitive polynomial over the integers modulo 2, or
in other words, a polynomial whose coefficients are either 0 or 1.
Suppose $P$ is such a polynomial of degree $q$:
\begin{equation}\label{6.4.2}
P=x^{q}+a_{1}x^{q-1}+a_{2}x^{q-2}+\dots +a_{q-1}x+1.
\end{equation}%
\noindent Define a sequence of integers $M_{k}$, by the $q$th term
recurrence relation:
\begin{equation}\label{6.4.3}
M_{k}=2a_{1}M_{k-1}\oplus 2^{2}a_{2}M_{k-2}\oplus \dots \oplus
2^{q-1}M_{k-q+1}a_{q-1}\oplus \left( 2^{q}M_{k-q}\oplus M_{k-q}\right) .
\end{equation}%
\noindent Here $\oplus $ denotes the XOR operation. The starting values for
the recurrence are $M_{1},\dots ,M_{q}$ that are odd integers chosen
arbitrarily and less than $2,\dots ,2^{q}$, respectively. The directional
numbers $V_{k}$ are given by:
\begin{equation}\label{6.4.4}
V_{k}=\frac{M_{k}}{2^{k}}\qquad k=1,\dots ,w.
\end{equation}%
\noindent Table \ref{DirNumb_S} shows the first ten primitive
polynomials and the starting values used to generate the direction
numbers for the $10$ dimensional Sobol' sequence.
\begin{figure}
\begin{center}
\includegraphics[width=0.7\textwidth]{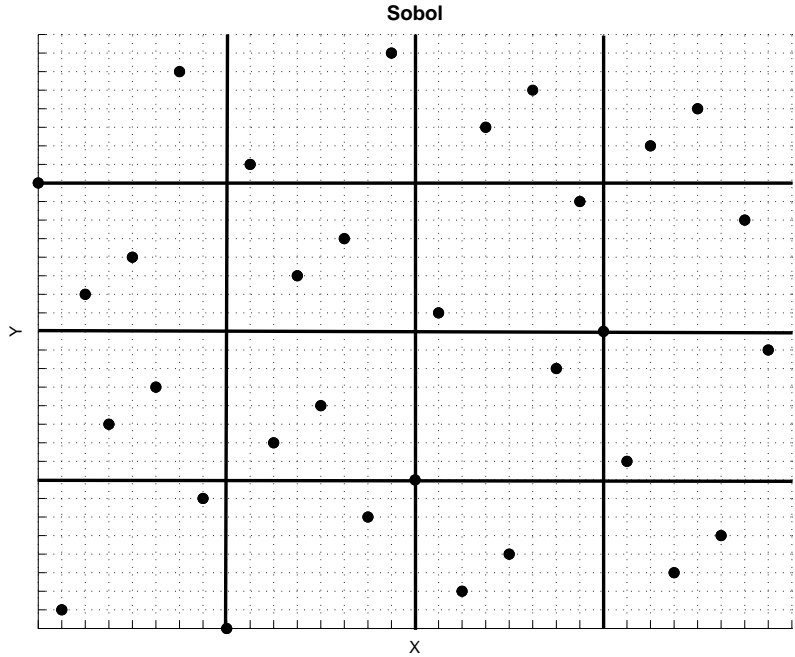}
\caption{The panel shows the first 32 points of the 2-dimensional Sobol' sequence}
\end{center}
\end{figure}

%
%
\subsection{Scrambling Techniques}

\noindent Digital nets are deterministic sequences. Their properties
ensure good distribution in the hypercube $\left[ 0,1\right) ^{d}$,
enabling precise sampling of all random variables, even if they are
very skewed. The main problem is the computation of the error in the
estimation, since it is difficult to compute and depends on the
chosen integrand function. To review, the crude MC provides an
estimation with low convergence independent of $d$ and the
possibility to statistically evaluate the RMSE. On the other hand,
the QMC method gives a higher convergence, but there is no way to
statistically calculate the error.

In order to estimate a statical measure of the error of the
Quasi-Monte Carlo method we need to randomize a $\left( t,m,d\right)
$-net and try to
obtain a new version of points such that it still is a $\left( t,m,d\right) $%
-net and has uniform distribution in $\left[ 0,1\right) ^{d}$.

This randomizing procedure is called scrambling. The scrambling technique
permutes the digits of the digital sequence and returns a new sequence that
has both the properties described above.

The scrambling technique we use is called Faure-Tezuka Scrambling
(for a precise description see Owen \cite{ow2002}, Hong and
Hickernell \cite{HH2000}).

For any $z\in \left[ 0,1\right) $ we define $\mathbf{\Psi }(z)$ as the $%
\infty \times 1$ vector of the digits of $z$.

Now let $L_{1},\dots L_{d}$ be nonsingular lower
triangular $\infty \times \infty $ matrices and let $\mathbf{e_{1}},\dots ,%
\mathbf{e_{d}}$ be $\infty \times 1 $ vectors. Only the diagonal elements of
$L_{1},\dots L_{d} $ are chosen randomly and uniformly in $%
Z_{b}^{\ast }=\left\{1,\dots,b\right\}$ , while the other elements
are chosen in $Z_{b}=\left\{0,1,\dots,b\right\}$. $\mathbf{Y}$, the
Faure-Tezuka scrambling version of $\mathbf{X}$, is defined as:

\begin{equation}\label{6.5.1}
\mathbf{\Psi }\left( y_{ij}\right) =\left( L_{j}\,\mathbf{\Psi }%
\left( x_{ij}\right) +\mathbf{e}_{j}\right) mod b
\end{equation}
\noindent All operations take place in the finite field $Z_{b}$.
Owen proved that, with his scrambling, it is possible to obtain
(see Owen \cite{ow2003}):
\begin{equation}\label{6.5.2}
Var\left[ \hat{I}\right] \leq \frac{b^{t}}{n}\left[ \frac{b+1}{b-1}\right)
^{d}\sigma ^{2},
\end{equation}%
\noindent for any twice integrable function in $\left[ 0,1\right) ^{d}$. These
results state that for low dimension $d$ , the randomized QMC (RQMC)
provides a better estimation with respect to Monte Carlo, at least for large $%
n $.
\begin{figure}
\begin{center}
\includegraphics[width=0.7\textwidth]{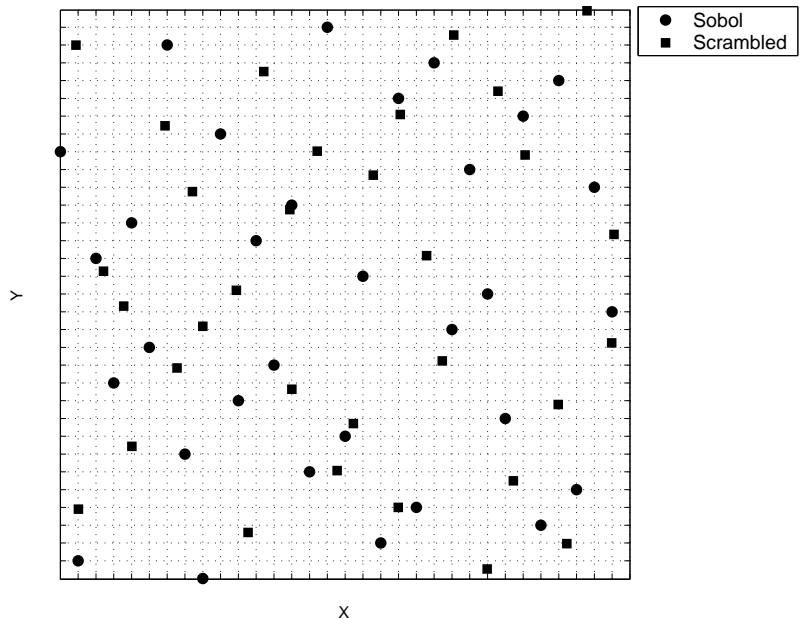}
\caption{The panel shows the first 32 point of the Sobol sequence
compared to their Faure-Tezuka scrambled version}
\end{center}
\end{figure}
%
%
\section{Implementation and Algorithm}
\noindent We illustrate the simulation procedure to compute the
arithmetic Asian option price. The purpose of our analysis is to
characterize the efficiency of Monte Carlo methods based on the path
generation techniques and the uniform points used for the evaluation
of the integral (\ref{2.2.10}). We consider separately the constant
volatility and time-dependent volatility markets.

It must be stressed that Quasi-Monte Carlo estimations
are dramatically influenced by the problem dimension, because the rate of
convergence depends on the problem dimension $d$, as it can be seen in equations (%
\ref{6.2.1}) and (\ref{6.5.2}). Many studies and experiments suggest
that Quasi-Monte Carlo methods can only be used for problem
dimensions up to $20$ (see Boyle, Broadie and Glasserman
\cite{BBG1995} for more on this issue). This condition translates
into a relationship between the number $M$ of underlying assets and
the number $N$ of monitoring times: $M\times N\leq 20$. When this
condition is not satisfied anymore we use the Latin Supercube method
that we describe hereafter.
%
%
\subsection{Latin Supercube Sampling}
The scrambling procedure allows the statistical estimation of the
RMSE as the crude MC does with the order of convergence that depends
on the the dimension $d$. For high $d$ the fast convergence of the
RQMC is lost, there is no benefit to use it compared to the simple
MC. Generally in finance the dimension is high even using dimension
reduction techniques like ANOVA-PCA decomposition.
\\ Owen \cite{ow1998B} has proposed a method to extend the convenience of applicability of
RQMC for high dimensions. This method is called \emph{Latin
Supercube Sampling}, (LSS), owing to its similarity to the LHS. The
random permutation is now applied to a set of subsequence of the
original one with some statistical sense.
\\ Let $\mathbf{Y}=\{\mathbf{y_1},\dots,\mathbf{y_{b^{m}}\}}$ be the
digital sequence of the simulation variables, and ${b^{m}}=N$.
Dividing it into $k$ nonempty and disjoint subsets
$\mathbf{Y}=\bigcup_{r=1}^k \mathbf{Y_r}$ and letting
$s_r=dim\mathbf{Y_r}$ we have $\sum_{r=1}^k s_r=d$. In practice,
 each point of the sequence can be represented as
 $\mathbf{y_i}=(\chi_i^1,\dots,\chi_i^k)$, where
 $\chi_i^r\in[0,1[^{s_r}$; these points $\chi_i^r$ are ordinarily
 points of an $s_r$-dimensional RQMC method.
 \\ For $r=1,\dots,k$ let $\pi_r(i)$ be an independent uniform and random permutation of $\{1,\dots,N\}$
 than a Latin Supercube sample is obtained by taking:
 \begin{equation}\label{7.1}
 \mathbf{\hat y_i} =(\chi_{\pi_1(i)}^1,\dots,\chi_{\pi_k(i)}^k)
\end{equation}
It means that the first $s_1$ colums in the LSS are
obtained by randomly permuting the run order of the RQMC
points $\chi_i^1,\dots,\chi_i^1$, the next $s_2$ columns come from
an independent permutation of the run order of $\chi_i^2$ and so
on.

The convenient way to divide the original set might be arranging
them in statistically orthogonal sets using the ANOVA-PCA
decomposition.

 In practice in financial simulation with $d$
Brownian motions, it may make sense to select $5$ principal
components of each path, to apply an RQMC method to each of them
with LSS and then pad them out the other variables with LHS. In
fact, about  $95\%$ of the total variance of the Brownian motion is
explained by these components. Alternatively, it may be better to
group the first $k$ principal component, then the second, and so on.

However, all these results are weak and only the practical test can
give an answer to which sequence and scrambling should be used.
%
%
\subsection{Key Steps of the Simulation Procedure}
As a first scenario we run simulations using the Cholesky and the
PCA decomposition procedures for the constant volatility case. As a
second scenario we test the efficiency of the proposed Kronecker
product approximation by comparing the its results with those
obtained with the PCA decomposition.

As a random number generator we use three configurations: standard,
LHS and Faure-Tezuka scrambled version of the Sobol' sequence.

The test for constant volatility consists of three main steps:
\begin{enumerate}
  \item Random number generation by standard MC, LHS or RQMC.
  \item Path generation with Cholesky and PCA decompositions.
  \item Monte Carlo estimation.
\end{enumerate}
For the time-dependent volatility case the three steps are:
\begin{enumerate}
  \item Random number generation by RQMC.
  \item Path generation with PCA decomposition and KPA.
  \item Monte Carlo estimation.
\end{enumerate}
The first step of both cases is realized by using the correspondent
random generator of uniform random variables. In order to extract
normal random variables we rely on the inverse transform method that
require the numerical inversion of the cumulative function of the
standard normal. This numerical procedure may destroy the better
stratification and the uniformity introduced by LHS and especially
by low-discrepancy sequences. We use the Moro's algorithm that is
more precise than the standard one due to Beasley and Springer. It
provides a better accuracy on the tails of the inverse normal where
we require that the LHS and Sobol sequences must reveal their higher
precision, see Moro \cite{Moro1995} and Glasserman \cite{Glass2004}
for more on the topic.

For constant volatilities the second step can be implemented by the
following algorithm:
\begin{enumerate}
\item Define the parameters of the simulation.
\item Define the drift as in equation (\ref{2.2.7}).
\item Create the $N\times N$ correlation matrix $\left(R%
\right)_{l,k}=\left(t_{l}\wedge t_{k}\right);l,k=1,\dots ,N$.
\item Define the correlation matrix $\Sigma$ among the $M$ asset returns.
\item Perform either a PCA or the Cholesky decomposition on the
global correlation matrix $\Sigma_{MN}$. This matrix is built up
by repeating the constant block of correlation $\Sigma$ at all the
times of observation.
\end{enumerate}
For time-depending volatilities we define the drift as equation
(\ref{2.2.8}), while the last operation consists of performing the
PCA decomposition and the KPA.

Stratification introduces a correlation among random drawings so
that the hypothesis of the Central Limit theorem is not satisfied
and we cannot compute the RMSE straightforward.
We rely on the batch method that consists of repeating $%
N_{B}$ simulations for $B$ times (batches).
We assume that each of the $B$ batches eliminates the correlation and the results form
a sequence of $B$ independent random variables. We compute the average Asian
price for each batch; the RMSE becomes:
\begin{equation}\label{7.2.1}
RMSE=\sqrt{\frac{\sum_{b=1}^{B}\left( \bar{a}\left( 0\right) _{b}-\bar{a}%
\left( 0\right) \right) ^{2}}{B\left( B-1\right) }},
\end{equation}%
\noindent where $\left( \bar{a}\left( 0\right) _{1},\dots ,\bar{a}\left(
0\right) _{B}\right) $ is a sample of the average present values of the
Asian option generated in each batch.
%
%
\section{Numerical Experiments}

\noindent We perform a test of all the valuation procedures
described in the previous section. We specify our investigation into
constant and time-dependent volatilities cases while our experiments
involve standard Monte Carlo, the Latin Hypercube Sampling and
Randomized Quasi Monte Carlo with the Faure-Tezuka scrambled version
of the Sobol' sequences.
%
%
\subsection{Constant Volatility: Results}
As for a first pricing experiment we consider an at-the-money
arithmetic Asian option with strike price $K=100$, written on a
basket of $M=2$ underlying assets, expiring at $T=1$ year and
sampled $N=5$ times during its lifetime.

All results are obtained by using $%
S=8192$  drawings and $10$ replications. Table \ref{test_S}
reports the input parameters for our test.
\begin{table}[tbp] \centering%
\begin{tabular}{||c||}
\hline\hline
$%
\begin{array}{lll}
S_{i}\left( 0\right) & = & 100 \\
K & = & 100 \\
r & = & 2\% \\
T & = & 1 \\
\sigma _{1} & = & 30\%\quad \\
\sigma _{2} & = & 40\% \\
\rho _{ij} & = & 0\text{ and }40\%\quad \text{for }i,j=1,2.%
\end{array}%
$ \\ \hline\hline
\end{tabular}%
\caption{Input Parameters Used in the First Simulation}\label{test_S}%
\end{table}
 The nominal dimension of the problem is $M\times N=10$ that is
equal to the number of rows and columns of the global correlation matrix $%
\Sigma_{MN}$. Paths are simulated by using both PCA and the
Cholesky decomposition as in Dahl and Benth \cite{DB2001} and
\cite{DB2002}.

Table \ref{UnCorr_S} and Table \ref{Corr_S} show the results for the
positive correlation and uncorrelated cases, respectively. Simulated
prices of the Asian basket options are in statistical accordance,
while the estimated RMSEs depend on the sampling strategy adopted.
The rate of convergence of the RQMC estimation is higher than the
other two methods. In particular it is ten times higher than the
standard Monte Carlo method that would return the same accuracy with
$100\times S$ drawings.

We observe that the PCA generation provides a better estimation both
for LHS and RQMC, because these ones are more sensitive to the
effective dimension, while PCA causes no distinction for the
standard MC. The effect is more pronounced for the correlation case
where the more complex structure of the global correlation matrix
$\Sigma_{MN}$ influences the estimation procedure.

As from a financial perspective, it is normal to find a higher price
in the positive correlation case than in the uncorrelated one.

\begin{table}[tbp] \centering
\begin{tabular}{||l|l|l|l||}
\hline\hline
& \textbf{Standard MC} & \textbf{LHS} & \textbf{RQMC} \\ \hline
\textbf{PCA} & 7.195 (0.016) & 7.157 (0.013) & 7.1696 (0.0017) \\
\textbf{Cholesky} & 7.242 (0.047) & 7.179 (0.022) & 7.1689 (0.0071) \\
\hline\hline
\end{tabular}%
\caption{Uncorrelation Case. Estimated Prices and Standard Errors.}\label{UnCorr_S}%
\end{table}%

\begin{table}[tbp] \centering%
\begin{tabular}{||l|l|l|l||}
\hline\hline
& \textbf{Standard MC} & \textbf{LHS} & \textbf{RQMC} \\ \hline
\textbf{PCA} & 8.291 (0.053) & 8.2868 (0.0073) & 8.2831 (0.0016) \\
\textbf{Cholesky} & 8.374 (0.055) & 8.293 (0.026) & 8.2807 (0.0064) \\
\hline\hline
\end{tabular}%
\caption{Correlation Case. Estimated Prices and Standard Errors.}\label%
{Corr_S}%
\end{table}
Moreover, we develop our analysis by investigating a very
high-dimensional pricing problem. A basket of $M=10$ underlying
assets is considered with $N=250$ sampling time points, the nominal
dimension is $d=2500$.

We run our simulation with the same parameters used by Imai and Tan
\cite{IT2002} and  use the LSS for the high-dimensional QMC
estimation as presented in the cited reference. The authors
concatenated 100 or 50 sets of 25 or 50 dimensional Sobol' sequence,
respectively. They exploit the LSS method in order to obtain a
complete 2500 dimensional sample of digital net. Owen \cite{ow1998}
is more restrictive; the author suggests to use scrambled digital
sequences for the first five or ten components and LHS for the
others or to concatenate the principal components. We compare the
results and investigate the effective dimensions and the
contribution of the eigenvalues of the global correlation matrix.
Table \ref{test_S2} reports input parameters for our test.
\begin{table}[tbp] \centering%
\begin{tabular}{||c||}
\hline\hline
$%
\begin{array}{lll}
S_{i}\left( 0\right) & = & 100 \\
K & = & 100 \\
r & = & 4\% \\
T & = & 1 \\
\sigma_i & = & 10\%+\frac{i-1}{9}40\%\quad \text{for }i=1,\dots,10
\\\rho_{ij} & = & 0\quad \textrm{and}\quad40\%\quad \text{for }i,j=1,\dots,10
\end{array}%
$ \\ \hline\hline
\end{tabular}%
\caption{Input Parameters Used in the Second Simulation}\label{test_S2}%
\end{table}

We compute the eigenvalues and eigenvectors of $\Sigma_{MN}$.
Property (10) of the Kronecker product is fundamental in this
computation and considerably reduces the computational burden and
time. It results that the effective dimension is $143$ or $170$ for
the correlation and uncorrelation cases, respectively, are much
smaller than the nominal one. Considering the first $143$($170$)
columns, that is the first $143$($170$) principal components, the
generating matrix $C$ takes into account $99\%$ of the total
variance.

Table \ref{test} shows all the results we obtained.  We concatenate
$50$ sets of $50$-dimensional randomized low-discrepancy sequences.
%
%

\begin{table}[tbp] \centering
\begin{tabular}{||l|l|l|l||}
\hline\hline \textbf{Uncorrelation} & \textbf{Standard MC} &
\textbf{LHS} & \textbf{RQMC}
\\ \hline
\textbf{PCA} &  3.414(0.015) &  3.4546(0.0054) & 3.4438(0.0015) \\
\textbf{Cholesky} & 3.426(0.015) &  3.4323(0.0070) &  3.4518(0.0058) \\
\hline\hline  \textbf{Correlation}& \textbf{Standard MC} &
\textbf{LHS} & \textbf{RQMC}
\\ \hline
\textbf{PCA} &  5.648(0.029) &  5.6655(0.0032) &  5.65750(0.00040) \\
\textbf{Cholesky} &  5.604(0.029) &  5.670(0.013) &  5.63710(0.019) \\
\hline\hline
\end{tabular}%
\caption{Prices and RMSEs both for the correlated and uncorrelated case when $100\%$ of the variance is considered.}\label{test}%
\end{table}
We  consider both the matrix of $2500$ rows and $143$($170$)
columns, excluding the effects of the remaining principal
components, and the complete ANOVA in order to investigate the
effectiveness of our assumptions and hypotheses.

Table \ref{eig} presents the different Monte Carlo estimations with
respect to the number of eigenvalues when LHS is used.
\begin{table}[tbp] \centering%
\begin{tabular}{||c|c|c|c|c|c||} \hline\hline
    {\bf } & {\bf Positive Correlation} &     {\bf } &     {\bf } & {\bf Zero Correlation} &     {\bf } \\
\hline {\bf Price} & {\bf RMSE} &    {\bf E} & {\bf Price} & {\bf
RMSE} &    {\bf E} \\ \hline
   5.262 &    0.090 &          5 &    2.596 &    0.041 &          5 \\
   5.294 &    0.088 &         10 &    3.190 &    0.047 &         10 \\
   5.433 &    0.088 &         15 &    3.212 &    0.047 &         15 \\
   5.528 &    0.091 &         20 &    3.239 &    0.047 &         20 \\
   5.484 &    0.092 &         25 &    3.289 &    0.047 &         25 \\
   5.445 &    0.090 &         30 &    3.375 &    0.048 &         30 \\
    5.653 &    0.015 &         147 &    3.452 &    0.010 &       170 \\
\hline\hline
\end{tabular}
\caption{Prices and RMSE for different principal components when LHS
is used.}\label{eig}
\end{table}
\begin{table}[tbp] \centering
\begin{tabular}{||l|l|l|l||}
\hline\hline \textbf{Uncorrelation} & \textbf{RQMC} &
\textbf{Correlation} & \textbf{RQMC}
\\ \hline
\textbf{PCA} &  3.4475(0.0023) &  \textbf{PCA} & 5.65860(0.00072) \\
\textbf{Cholesky} & 3.426(0.0087) &  \textbf{Cholesky} &  5.603(0.022) \\
\textbf{LT} & 3.4461(0.0012) &  \textbf{LT} &  5.6780(0.00047) \\
\hline\hline
\end{tabular}%
\caption{Estimated Results by Imai and Tan
\cite{IT2007}\label{ImaiTan}}
\end{table}

Table \ref{ImaiTan} illustrates the values found by Imai and Tan
\cite{IT2002}. Their results were obtained assigning the importance
of each component (not anymore PCA) with their LT method. All the
estimations found are unbiased and in agreement with those presented
in the cited references.

The Quasi-Monte Carlo method with LSS extension proves to be a
powerful variance reduction technique, particularly when coupled
with the ANOVA-PCA decomposition. Moreover, the Kronecker product
turns out to be a fast tool to generate multi-dimensional Brownian
paths. Indeed, the elapsed time to realize the same path without
using the properties of the Kronecker product is a lot higher.

The estimation with Cholesky decomposition gives higher uncertainty
than the PCA approach, meaning that a small amount of variance is
lost. This is due to the fact that a relevant part of the variance
is carried out by a few eigenvalues of the  covariance matrix $R$.
If these eigenvalues are observed, it can be noticed that only few
of them are relevant in the PCA analysis and they are much bigger
than the ones of the matrix $\Sigma$.
\subsection{Constant Volatility: Comments}
Based on these results, we can make the following conclusions:
\begin{enumerate}
  \item The RQMC method and the use of the Faure-Tezuka scrambling
technique provide the best estimation among all the implemented
procedures for both the \textquotedblleft
Correlation\textquotedblright\ and \textquotedblleft Zero
Correlation\textquotedblright\ cases. The correspondent RMSEs are
the smallest ones with a higher order of convergence with the same
number of simulations.
  \item The Kronecker product is a fast and efficient tool for
generating multi-dimensional Brownian paths with a low computational
effort.
  \item As compared to to the standard Monte Carlo and LHS approaches, the
use of scrambled low-discrepancy sequences provides more accurate
results, at least for $M\times N\leq 20$, particularly with the PCA
and LT-based methods.
  \item The accuracy of the estimates is strongly dependent on the
choice of the Cholesky or the PCA approach. In particular,
independent of the simulation procedure (MC, LHS or RQMC), when
using PCA decomposition the estimates are affected by a smaller
sampling error (smaller standard error).
\end{enumerate}

%
\subsection{Time-dependent Volatility: Results}
The constant volatility hypothesis is the starting point for the
pricing problem. A further improvement can be achieved by
considering a time-dependent volatility function.

It is market practice to choose step-wise time-dependent
volatilities. We want to investigate a more complex dependence to
test our new approach based on the Kronecker product approximation.
For this aim, we adopt an exponentially decaying function having the
following expression:
\begin{equation}
\sigma_i =
\hat{\sigma}_i(0)\,exp\big(-\frac{t}{\tau_i}\big)+\sigma_i(+\infty)
\end{equation}
where $\hat{\sigma}_i(0)+\sigma_i(+\infty)=\sigma_i(0)$ is the
initial volatility for the $i$-th asset, $\sigma_i(+\infty)$ is its
asymptotic volatility and $\tau_i$ its decay constant.

The particular time-dependent function leads to the following
solution:
\begin{eqnarray*}
\int_0^{t_j\wedge t_l} \sigma_i(t) \sigma_k(t)\rho_{ik}dt&=&
\hat{\sigma}_i(0)\hat{\sigma}_k(0)\tau_
{ik}\Big(1-exp\big(-\frac{t}{\tau_{ik}}\big)\Big)+\\ &&
  +\hat{\sigma}_i(0)\sigma_k(+\infty)\tau_{i}\Big(1-exp\big(-\frac{t}{\tau_{i}}\big)\Big)+\\
&& +
\hat{\sigma}_k(0)\sigma_i(+\infty)\tau_{ik}\Big(1-exp\big(-\frac{t}{\tau_{ik}}\big)\Big)+\\
&& +\sigma_i(+\infty)\sigma_k(+\infty)t
\end{eqnarray*}
where $\tau_{ik}=\tau_i\tau_k/(\tau_i+\tau_k)$.

The simulation implemented to obtain the  price of an Asian option
supposing time-dependent volatility evolves as the constant
volatility case. The main difference is the procedure to reduce the
dimension of the problem.

The parameters chosen for the simulation are listed in table
\ref{KPA1}.
\begin{table}[tbp] \centering%
\begin{tabular}{||c||}
\hline\hline
$%
\begin{array}{lll}
S_i(0) & = & 100
\\ r & = & 4\%
\\ T & = & 1 \textrm{year}
\\ \sigma_i(0) & = & 10\%+\frac{i-1}{9}40\%
\\ \sigma_i(+\infty) & = & 9\% \qquad\textrm{for all i}
\\ \tau_i & = & 1.5 \textrm{year}
\\ K & = & 100
\\ \rho_{ij} & = & 0\quad \textrm{and}\quad40\%\quad \textrm{for}\quad
i,j=1,\dots,10
\end{array}%
$ \\ \hline\hline
\end{tabular}%
\caption{Input Parameters for the Time-depending Case}\label{KPA1}
\end{table}

The initial volatilities are equal to those  used in the constant
volatility case. The asymptotic volatility and the decay constant
are the same among all the assets. These parameters are chosen in
order to allow a comparison with respect to the constant volatility
case. Indeed, the price of the options is sensitive to the change of
volatility and in particular its decreasing trend should provide a
lower price.

The basket consists of $10$ underlying assets, the time grid has
$250$ equally spaced points and the number of runs is $S=8192$ and
$10$ replications. Table \ref{KPA2} shows the results coming from
the simulation using the RQMC method both with the KPA and the PCA
for dimension reduction.
\begin{table}\centering
\begin{tabular}{||cc|cc||}
 \hline\hline  {\bf Positive Correlation (KPA)}  &            & {\bf Zero
Correlation(KPA)} &
\\
\hline
     {\bf Price} &   5.19658 &      {\bf Price} &   3.20784 \\
      {\bf RMSE} &   0.00063 &       {\bf RMSE} &   0.00040 \\
         {\bf E} &        145 &          {\bf E}  &        173 \\
\hline\hline {\bf Positive Correlation (PCA)}  &            & {\bf
Zero Correlation(PCA)} & \\
 \hline
     {\bf Price} &   5.19856 &      {\bf Price} &   3.20147 \\
     {\bf RMSE} &   0.00062 &       {\bf RMSE} &   0.00040 \\
     {\bf E} &        123 &          {\bf E}  &        150 \\
\hline\hline
\end{tabular}
\caption{Estimated Results for the Time-depending Case, ANOVA =
0.99}\label{KPA2}
\end{table}

The  KPA path-generation is efficient and fast. To have an idea of
its speed, the elapsed times to obtain the generating matrix without
exploiting the properties of the Kronecker product and no
approximations are more than ten times higher. As expected, the
simulation gives smaller prices with respect to the constant
volatility situation, because a decreasing volatility function has
been assigned.

The nominal dimensions of the problem $E$ using PCA come out to be
$126$ and $150$ for the correlation and uncorrelation cases. When
adopting the KPA the approximated nominal dimensions are higher,
 $145$ and $173$, respectively. If we would consider ANOVA $=0.9885$
for the correlation case and ANOVA = $0.98805$ for the uncorrelation
case we would get the PCA-found nominal dimension for ANOVA=$0.99$.
We can judge this small difference as negligible and consequently
our approximating technique to be efficient and leading to
consistent results. As with $N$ and $M$ small, the Cholesky
decomposition alone would require a small number of operations
without giving any order of the importance for the random sources.

\begin{table}[tbp] \centering
\begin{tabular}{||l|l|l|l||}
\hline\hline & \textbf{KPA} & \textbf{PCA} & \textbf{Cholesky}
\\ \hline
\textbf{Price} & 3.20545 & 3.20390 & 3.1838  \\
\textbf{RMSE} & 0.00040 & 0.00041 & 0.0091  \\
\hline\hline
\end{tabular}%
\caption{Uncorrelation Case. Estimated Prices and Standard Errors. ANOVA$=1$.}\label{Un_det}%
\end{table}%

\begin{table}[tbp] \centering%
\begin{tabular}{||l|l|l|l||}
\hline\hline & \textbf{KPA} & \textbf{PCA} & \textbf{Cholesky}
\\ \hline
\textbf{Price} & 5.20060 & 5.20210 & 5.1946 \\
\textbf{RMSE} & 0.00050 & 0.00058 & 0.0093 \\
\hline\hline
\end{tabular}%
\caption{Correlation Case. Estimated Prices and Standard Errors.
ANOVA$=1$.}\label{cor_det}
\end{table}

Tables \ref{Un_det} and \ref{cor_det} present the estimated prices
when taking into account the full components. All the results are in
accordance with those ones found with ANOVA$=0.99$.

Table \ref{KPA3} illustrates the sensitivity with respect to the
number of principal components $E$:
\begin{table}\centering
\begin{tabular}{||c|c|c|c|c|c||}
\hline\hline
    {\bf } & {\bf Positive Correlation} &     {\bf } &     {\bf } & {\bf Zero Correlation} &     {\bf } \\
\hline {\bf Price} & {\bf RMSE} &    {\bf E} & {\bf Price} & {\bf
RMSE} &    {\bf E} \\ \hline
   5.7805 &    0.0079 &          5 &    2.6368 &    0.0038 &          5 \\
   4.9904 &    0.0081 &         10 &    2.9681 &    0.0042 &         10 \\
   5.0226 &    0.0081 &         15 &    3.1172 &    0.0043 &         15 \\
   5.1103 &    0.0081 &         20 &    3.0979 &    0.0043 &         20 \\
   5.1826 &    0.0083 &         25 &    3.1051 &    0.0043 &         25 \\
   5.1937 &    0.0082 &         30 &    3.1514 &    0.0043 &         30 \\
\hline\hline
\end{tabular}
\caption{Prices and RMSEs for different principal components. Case:
RQMC}\label{KPA3}
\end{table}
As in the constant volatility case it can be seen that the
estimation is convergent.

Moreover, we launch a new simulation with the LHS technique with the
same set of parameters. We list
 the estimated results in table \ref{KPA4}.
\begin{table}\centering
\begin{tabular}{||c|c|c|c|c|c||}
\hline\hline
    {\bf } & {\bf Positive Correlation} &     {\bf } &     {\bf } & {\bf Zero Correlation} &     {\bf } \\
\hline {\bf Price} & {\bf RMSE} &    {\bf E} & {\bf Price} & {\bf
RMSE} &    {\bf E} \\ \hline
   4.874 &    0.016 &          5 &    3.121 &    0.089 &          5 \\
   5.093 &    0.016 &         10 &    3.118 &    0.085 &         10 \\
   5.097 &    0.016 &         15 &    3.122 &    0.086 &         15 \\
   5.131 &    0.016 &         20 &    3.072 &    0.086 &         20 \\
   5.145 &    0.016 &         25 &    3.163 &    0.088 &         25 \\
   5.201 &    0.016 &         30 &    3.110 &    0.089 &         30 \\
\hline\hline
\end{tabular}
\caption{Prices and RMSEs for different principal components. Case:
LHS}\label{KPA4}
\end{table}
The estimated prices have higher RMSEs, confirming the fact that the
RQMC approach provides a good variance reduction.
\subsection{Time-Dependent Volatility: Comments}
According to the results we have found in the time-dependent case,
it is possible to draw the following conclusions:
\begin{enumerate}
  \item RQMC with LSS is a general approach that does not depend on
  the chosen price dynamic.
  \item The KPA we propose, provides unbiased estimations with a
  reduction of the computational cost. In the framework we investigate,
  KPA returns a higher nominal dimension, as expected, but only relatively to
  a negligible amount of variance.
  \item KPA is a lot faster than the straightforward PCA because it
  exploits the properties of the Kronecker product and the \emph{boomerang}
  shaped matrices. The \emph{ad hoc} Cholesky decomposition
  algorithm we develop is fundamental for the KPA.
  We do not report computational times because we expect that
  further improvements can be done.
  \item KPA and PCA can be considered both valid as path-generation
  methods to support the ANOVA and the identifications of effective
  dimensions.
\end{enumerate}

\end{document}